\documentclass[leqno,11pt]{article}
\usepackage{amsmath, amsthm}
\usepackage{amsfonts}
\usepackage{amssymb}
\usepackage{graphicx}
\usepackage[monochrome]{color}
\usepackage{enumitem}

\newtheorem{theorem}{Theorem}[section]

\newcommand{\RR}{\mathbb{R}}
\newcommand{\p}{{\delta^\gamma}}

\newcommand{\re}{\mathbb{R}}
\newcommand{\ren}{\mathbb{R}^N}
\newcommand{\ve}{\varepsilon}

\newcommand{\dx}{\,{\rm d}x}
\newcommand{\dy}{\,{\rm d}y}

\newcommand{\rd}{{\rm d}}
\newcommand{\A}{\mathcal{L}}
\def\dist{\mathrm{dist}}
\numberwithin{equation}{section}
\def\qed{\,\unskip\kern 6pt \penalty 500
\raise -2pt\hbox{\vrule \vbox to8pt{\hrule width 6pt
\vfill\hrule}\vrule}\par}

\title{\bf The mathematical theories of diffusion.\\
Nonlinear and fractional diffusion\footnote{To appear in Springer Lecture Notes in Mathematics, C.I.M.E. Subseries, 2017.}}
\author{ {\Large Juan Luis V\'azquez} \\[10pt]
Universidad Aut\'onoma de Madrid, Spain}

\date{ \ } 
\begin{document}
\maketitle

\begin{abstract}
We describe the mathematical theory of diffusion and heat transport with a view to including some of the main directions of recent research. The linear heat equation is the basic mathematical model that has been thoroughly studied in the last two centuries. It was followed by the theory of parabolic equations of different types. In a parallel development, the theory of stochastic differential equations gives a foundation to the probabilistic study of diffusion.

Nonlinear diffusion equations have played an important role not only in theory but also in physics and engineering, and we focus on a relevant aspect, the existence and propagation of free boundaries. Due to our research experience, we use the porous medium  and fast diffusion equations as  case examples.

A large part of the paper is devoted to diffusion driven by fractional Laplacian operators and other nonlocal integro-differential operators representing nonlocal, long-range diffusion effects. Three main models are examined (one linear, two nonlinear), and we report on recent progress in which the author is involved.
 \end{abstract}

\tableofcontents
\section{\textsf{Introduction to Diffusion}}
\label{sec.intro}


There are a number of phenomena in the physical sciences that we associate with the idea of diffusion. Thus, populations of different kinds diffuse; particles in a solvent and other substances diffuse. Besides, heat propagates according to a process that is mathematically similar, and this is a major topic in  applied science. We find many other instances of diffusion: electrons and ions diffuse; the momentum of a viscous  fluid diffuses (in a linear way, if we are dealing with a Newtonian fluid). More recently, we even talk about diffusion in the financial markets.

The word diffusion derives from the Latin  {\sl diffundere}, which means ``to spread out''. A substance spreads out by moving from an area of high concentration to an area of low concentration. This mixing behaviour does not require any bulk motion, a feature that separates  diffusion from other transport phenomena like convection, or advection.

\medskip

\begin{center}
 \includegraphics[width=7cm]{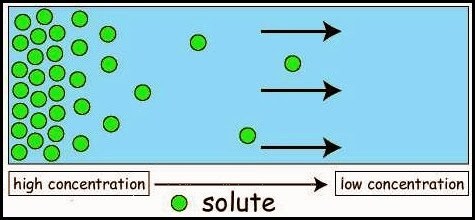}
\end{center}

 This description can be found in {\sl  Wikipedia} \cite{WikiDiff}, where we can also find a longer listing of more than 20 diffusive items,  that are then developed as separate subjects. It includes all classical topics we mention here.

\subsection{Diffusion in Mathematics}

In this survey paper we want to present different topics of current interest in the mathematical theory of diffusion in a historical context. To begin with, we may ask if mathematics is really relevant in the study of diffusion process? The answer is that diffusion is a topic that enjoys superb mathematical modelling. It is a branch of the natural sciences that is now firmly tied to a number of mathematical theories that explain its working mechanism in a quite successful way. The quantity that diffuses can be  a concentration, heat, momentum, information, ideas, a price,... every such process can be called a diffusion, and its evolution is governed by mathematical analysis.

Going into the details of how we actually explain diffusion  with mathematics, it so happens that we may do it in a twofold way: roughly speaking, by means of the diffusion equation and its relatives, or by a random walk model and its relatives.

The older work concerns the description of heat propagation and mass diffusion by means of partial differential equations (PDEs), and this is the view that we are going to favor here. The type of PDEs used is the so-called parabolic equations, a family based on the properties of the most classical model, the linear Heat Equation (HE), which is called in this context the diffusion equation.

On the other hand, according to probabilists diffusion is described by random walks, Brownian motion, and more generally,  by stochastic processes, and this is a long and successful story in 20th century mathematics, culminated by  It\=o's calculus. Let us recall that the connection between the two visions owes much to A. N. Kolmogorov.

Actually, an interesting question for the reader or the expert is `How much of the mathematics of diffusion can be explained with {\sl linear models,} how much is {\sl essentially nonlinear}?' Linear models have priority when applicable by virtue of their rich theory and easier computation. But nonlinear models are absolutely necessary in many real-world contexts and most of our personal research has been based on them. Diffusion equations involving nonlinearities and/or nonlocal operators representing long-range interactions are the subject matter of the recent research that we want to report in this paper.

\medskip

\noindent {\sc Outline of the paper.} The declared intention is to make a fair presentation of the main aspects of Mathematical Diffusion as seen by an expert in PDEs, and then to devote preference to the work done by the author and collaborators, especially the more recent work that deals with free boundaries and  with fractional operators. The main topics are therefore the heat equation, the linear parabolic theory and the fractional diffusion in a first block; the nonlinear models, with emphasis on those involving free boundaries come next; finally, the nonlocal and nonlinear models, and here we will focus on the two equations that have been most studied by the author in the last decade, both combine porous medium nonlinearities and fractional diffusion operators.

Lengthy details are not frequent, but we give some for very recent work of ours and our collaborators.  On the other hand, we supply many important explicit solutions and comment on their role. Indeed, many such examples belong to the class {\sl stable diffusive patterns}, that combine their surprising occurrence in numerous real-world applications with the beauty of pure mathematics.

A large number of connections with other topics is given in the text, as well as  hints for further reading. More detail on the topics is to be found in the articles, monographs, or in our previous survey papers.

\medskip

\noindent {\sc Disclaimer.} Let us comment on an important absence. The stationary states  of diffusion belong to an important world, the {\sl elliptic equations}. Elliptic equations, linear and nonlinear, appear in a large number of applications: diffusion,
fluid mechanics, waves of all types, quantum mechanics, ... Elliptic equations are mathematically based on the Laplacian  operator, $\Delta=\nabla^2$, the most important operator for our community. This is a huge world. We are not going to cover in any detail the many developments in elliptic equations related to diffusion in this paper, we will just indicate some important facts and connections here and there.


\subsection{Heat equation. Main model for diffusion}

We begin our presentation with the linear heat equation (HE):
$$
\fbox  {$  u_t=\Delta u$}
  $$
proposed by J. Fourier as a mathematical model for heat propagation (``Th\'eorie Analytique de la Chaleur'', 1822, \cite{Fou1822}, with a previous attempt in 1807), and the Fourier analysis that he promoted. For a long time the mathematical study of heat transport and diffusion was almost exclusively centered on the heat equation. In these  two centuries, the mathematical  models of  heat propagation and  diffusion  have made great progress both in theory and application. Actually, they have had a strong influence on no less than 6 areas of Mathematics:  PDEs, Functional Analysis, Infinite-Dimensional Dynamical Systems, Differential Geometry and Probability, as well as  Numerics. And the theory has been influenced by its motivation from Physics, and in turn the concepts and methods derived from it have strongly influenced Physics and Engineering. In more recent times this influence has spread further away, to Biology, Economics, and the Social Sciences.

\medskip

\noindent $\bullet$ The classical analysis of the heat flow is based on two main
mathematical techniques: {\sl integral  representation} (convolution
with a Gaussian kernel) and {\sl Fourier analysis},  based on mode separation, analysis, and synthesis. Since this topic is well-known to the readers, see for instance \cite{EvansPDE}, we will stress the points that interest us to fix some concepts and tools.

\noindent {\bf 1. The  heat equation semigroup and Gauss.}
 When heat propagates in {\sl free space} $\ren$, the natural problem is the initial value problem
\begin{equation}
u_t=\Delta u, \qquad u(x,0)= f(x)\,,
\end{equation}
which is solved by convolution with the evolution version of the Gaussian function
\begin{equation}\label{gauss.kernel}
G(x,t)=(4\pi t)^{-N/2}\mbox{\rm exp}\,(-|x|^2/4t).
\end{equation}
Note that $G$ has all nice analytical properties for $t>0$, but for $t=0$ we have  $G(x,0)=\delta(x)$, a Dirac mass. $G$ works as a {\sl kernel}, a mathematical idea that goes back to Green and Gauss.

\smallskip

  The maps $S_t: u_0\mapsto u(t):=u_0\ast G(\cdot,t)$ form  a {\sl linear continuous semigroup } of contractions
in all $L^p$ spaces for all $1\le p\le \infty$. This is pure Functional Analysis, a product of the 20th century.

\smallskip

{\bf Asymptotic  behaviour as $t\to\infty$, convergence to the Gaussian}. If $u_0$ is an integrable function and $M=\int u_0(x)\,dx\ne 0$,  the following convergence is proved
\begin{equation}
\lim_{t\to\infty} t^{N/2} (u(x,t)-M\,G(x,t))= 0\,,
\end{equation}
and the limit holds uniformly in the whole space.  For convergence in $L^p$ less is needed. So $u(\cdot, t) $ increasingly resembles (i.e. as $t$ grows to infinity) a multiple of the Gaussian profile $G(\cdot, t)$. This is the famous {\sl Central Limit Theorem} in its continuous form (famous in Probability with $M=1$, but $M\ne 1$ makes no difference as long as $M$ is not zero).

The Gaussian function is the most famous example of the many  diffusive patterns that we will encounter, and the previous theorem shows that is not only stable but also asymptotic attractor of the heat flow (for finite-mass data). Note that the sharp convergence result needs renormalization in the form of the growth factor $t^{N/2}$.

\medskip

\noindent $\bullet$ {These are two classical personalities of the diffusion equation. }

\begin{center}
\includegraphics[width=5cm]{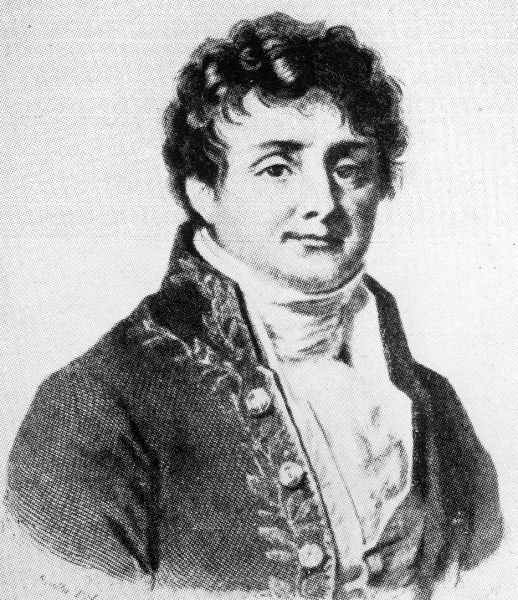}\qquad
\includegraphics[width=5cm]{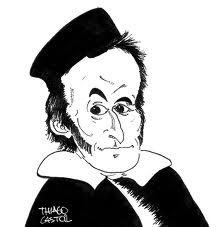}
\end{center}

\centerline{J. Fourier and K. F. Gauss}

\smallskip

{\bf 2. Matter diffusion.} This is an older subject in Physics, already treated by Robert Boyle in the 17th century with the study of diffusion in solids. After the work of Fourier in heat propagation,  Adolf Fick proposed his law of diffusion of matter \cite{Fick}, where the mass flux is proportional to the gradient of concentration and goes in the direction of lower concentrations. This leads to the heat equation, HE, as mathematical model. He also pointed out the fundamental analogy between diffusion, conduction of heat, and also electricity. Actually, Fourier's law for heat conduction (1822),  Ohm's law for electric current (1827),  Fick's law for diffusion in solids (1855), and Darcy's law for hydraulic flow (1856) have a similar mathematical gradient form.

\medskip

{\bf 3. The connection with Probability.} The time iteration of independent random variables with the same distribution led  to the theory of random walks, at the early times of the Bernoullis et al. Soon it was realized that this led to a large-time limit described (after renormalization) by the Gaussian distribution, as textbooks in Probability and Stochastic Processes show. The connection of this evolution with the heat equation took place after the construction of the {\bf Brownian motion} as a mathematically rigorous object in the form of {\sl Wiener process}. In the 1930s Kolmogorov investigated the equivalence of the two view points, i.e., the PDE approach via the Heat Equation and the stochastic approach via Brownian motion. This topic is covered by many PDE authors, let us mention \cite{EvansSt} and \cite{SalsaBk}.


Next we show two opposing diffusion graphs. They show the comparison of ordered dissipation in the heat equation view, as the spread of a temperature or concentration, versus the underlying chaos of the random walk particle approach, origin of the Brownian motion favoured by the probabilistic school.

\begin{center}
\includegraphics[width=5cm]{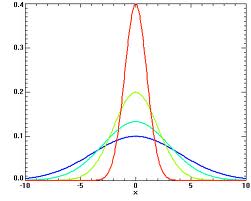}\qquad\qquad
\includegraphics[width=3.8cm]{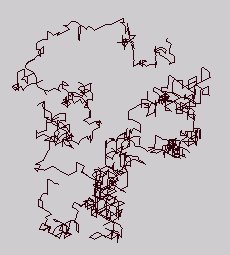}
\end{center}

 \centerline{\qquad Left, the nice HE evolution of a Gaussian
  \qquad Right, a sample of random walk\quad \quad}

\medskip

\noindent We will go back to the latter view in a while. The experimental observation of chaotic movement in Nature due to mechanical effects at the microscopic level is credited to Robert Brown (1827), see \cite{Brown}, hence the label `Brownian motion'.

\medskip

\noindent {\bf 4. The  Fourier Analysis approach to heat flows in bounded domains.} The second classical scenario for heat flows occurs when heat propagates inside a bounded domain of  space. The convolution approach does not work and other ideas have to be proposed. The Fourier approach proposes to look for a solution in the series form
\begin{equation}
u(x,t)=\sum T_i(t)X_i(x)\,,
\end{equation}
and then the time factors are easily seen to be negative exponentials of $t$, while the space components $X_i(x)$ form the spectral sequence, solutions of the problems
\begin{equation}
\fbox {$-\Delta X_i=\lambda_i \, X_i$}
\end{equation}
with corresponding eigenvalues $\lambda_i$. Boundary values are needed to identify the {\sl spectral sequence}  $(\lambda_i,X_i)$, $i=1,2,\dots$  This is the famous linear {\sl eigenvalue problem}, the starting point of the discipline of {\sl Spectral Theory.} The time-space coupling  implies then that $T_i(t)=e^{-\lambda_i t}$.  This is nowadays one of the most celebrated and useful topics in Applied Mathematics and is covered in all elementary PDE books.

The scheme works for many other equations of the form $u_t=A(u)$ and in this way  Fourier Analysis and Spectral Theory   developed  with great impetus, as well as  {Semigroup Theory.} Fourier analysis also took a direction towards the delicate study of functions, a proper branch of pure mathematics, which is one of the most brilliant developments in the last two centuries. Through the work of Cantor this also motivated advances in Set Theory, since the sets of points where a Fourier series does not converge can be quite complicated.

\medskip

\subsection{Linear heat flows}

We now consider a big step forward  in the mathematical tools of diffusion. A more general family of models was introduced to represent diffusive phenomena under less idealized circumstances and this was done both in the framework of PDEs and Probability. This happened in several stages.

In the framework of PDEs, the heat equation has motivated
the study of other linear equations, which now form the Theory of Linear Parabolic Equations. They are written in the form
\begin{equation}\label{linpar.eq}
\fbox  {$
u_t=\sum_{i,j}
a_{ij}\partial_i\partial_j u + \sum b_i \partial_i u + cu + f$}
\end{equation}
with variable coefficients $a_{ij}(x,t)$, $b_i(x,t)$, $c(x,t)$, and  forcing term $f(x,t)$. Belonging to the parabolic family requires some {\sl structure conditions} on such coefficients that will ensure that the solutions keep the basic properties of the heat equation. In practice, the main condition is the fact that the matrix $a_{ij}(x,t)$ has to be definite positive: there exists  $\lambda>0$
\begin{equation}\label{par.cond}
\sum_{i,j}a_{ij}(x,t)\xi_i\xi_j\ge \lambda \,\sum_i \xi_i^2\,.
\end{equation}
This must be valid for all vectors $\xi=(\xi_1,\dots,\xi_N)\in \ren$ and all $x,t$ in the space-time domain of the problem. Let us point out that the theory is developed under some additional conditions of regularity or size on the coefficients, a common feature of all PDE analysis. Another prominent feature is that all coefficients can be submitted to different more or less stringent conditions that allow to obtain more or less regular solutions.

A more stringent uniform positivity condition is
\begin{equation}\label{par.cond1}
\lambda \,\sum_i \xi_i^2 \le \sum_{i,j}a_{ij}(x,t)\xi_i\xi_j\le \Lambda \,\sum_i \xi_i^2
\end{equation}
with $0<\lambda<\Lambda $. Even the uniform condition \eqref{par.cond} can be relaxed so that  $\lambda>0$ depends on $x $ and $t$.  When these conditions are relaxed we talk about degenerate or singular parabolic equations. This flexibility on the structure conditions makes for a big theory that looks like an ocean of results. It will be important later in the nonlinear models.

In fact, the parabolic theory was developed in the sequel of its more famous stationary counterpart, the theory for the  elliptic equations,
\begin{equation}\label{linell.eq}
\sum_{i,j} a_{ij}\partial_i\partial_j u + \sum_i b_i \partial_i u + c\,u + f=0,
\end{equation}
with variable coefficients $a_{ij}(x), b_i(x),c(x)$ and forcing term $f(x)$.
The main structure condition is again \eqref{par.cond} or \eqref{par.cond1}, which is usually called uniform ellipticity condition. In the time dependent case, \eqref{par.cond} and \eqref{par.cond1} are called uniform parabolicity conditions.

Main steps in the {Parabolic Theory} are:

(1) The first step is the classical parabolic theory in which {\sl  $a_{ij}, b_i,c,f$} are assumed to be continuous or smooth, as needed. Functional spaces are needed as framework of the theory, and these turn out to be  $C^{\alpha}$ spaces
(H\"older) and the derived spaces  $C^{1,\alpha}$ and  $C^{2,\alpha}$. This leads to a well-known theory in which existence and uniqueness results, continuous dependence on data are obtained after adding initial and boundary data to the problem. And the theory provided us with  Maximum Principles, Schauder estimates, Harnack inequalities and other  very precise estimates.

This is a line of research that we would like to follow in all subsequent chapters when further models of diffusion will be treated, but unfortunately the direct application of the scheme will not work, and to be more precise, the functional setting will not be conserved and the techniques will change in a strong way.

The results are extended into the disciplines of Potential Theory and  Generation of Semigroups. These are also topics that will be pursued in the subsequent investigations.

\smallskip

(2) A first extension of the classical parabolic theory concerns the case where the {\sl coefficients are only continuous or bounded}. In the theory with bad coefficients there appears a bifurcation of the theory into Divergence and Non-Divergence Equations, which are developed with similar goals but quite different technology.
The difference concerns the way of writing the first term with second-order derivatives.

The way stated before is called non-divergent form, while the divergent form is
\begin{equation}\label{linell.div.eq}
\sum_{i,j} \partial_i(a_{ij}\partial_j u) + \sum b_i
\partial_i u + cu + f=0.
\end{equation}
This form appears naturally in many of the derivations from physical principles.
About the structure conditions, we assume the $a_{ij}$ to be bounded and satisfy the uniform parabolicity condition. The basic functional spaces are the Lebesgue and Sobolev classes, $L^p$,   $W^{1,p}$. Derivatives are understood as distributions or more often weak derivatives, and this motivates the label of {\sl weak theory}. Existence, uniqueness and estimates in $W^{1,p}$  or $W^{2,p}$ norms are produced. Maximum Principles, Harnack inequalities work and $C^{\alpha}$ is often proved. A very important feature is the Calder\'on-Zygmund theory, basic to establish regularity in Sobolev spaces.
Divergent form equations were much studied because of their appearance in problems of Science and Engineering. We refer to the books \cite{GT, LU} for the elliptic theory, and to \cite{Fr64, LSU, Lieb96} for the parabolic theory.

(3) There is nowadays a very flourishing theory of elliptic and parabolic equations with  bad coefficients in the non-divergence form \eqref{linpar.eq}, but we will not enter into it for reasons of space, since it does not affect the rest of our expos\'e.

\medskip

\subsection{The stochastic approach. SDEs}
 The probabilistic way to address the extension of the previous subsection  takes the form of the {\sl diffusion  process}, which  is a solution to a  stochastic differential equation, SDE for short. A diffusion is then a continuous-time Markov process with almost surely continuous sample paths. This is essentially a 20th century theory, originated in the work of Bachelier, Einstein, Smoluchowski, then Kolmogorov, Wiener and Levy, and the last crucial step was contributed by It\^o, Skorokhod, ... The stochastic equation reads
\begin{equation}\label{Ito.eq}
\fbox{ $dX=b\,dt+\sigma \,dW $}
\end{equation}
where $W$ is the $N$-dimensional Wiener process, and $b$ and $\sigma$ are (vector and matrix valued respectively) coefficients under suitable conditions. Differentials are understood in the It\^o sense. Among the extensive literature we mention Bass \cite{Bass}, Friedman \cite{Fr76}, Gihman-Skorohod \cite{GSk}, Oksendal \cite{Oksendal}, and Varadhan \cite{Varadhan} for the relation between PDEs and Stochastic processes. Thus, Bass discusses the solutions of linear  elliptic and parabolic problems by means of stochastic processes in Chapter II.  The solution of stochastic differential equations \eqref{Ito.eq} gives a formula to solve the Cauchy problem for the evolution PDE \ $u_t={\mathcal L}u$, where
$$
{\mathcal L}u = \sum a_{ij}\partial_{ij} u + \sum b_i \partial_i u,
$$
if   we take \normalcolor the vector  $b=(b_1,\dots, b_N)$ and the symmetric positive semi-definite matrix $a=(a_{ij})$ is given by   $a=\frac12 \sigma\cdot\sigma^T$. And there is an associated dual parabolic PDE in divergence form that is  solved by the same method. \normalcolor  See also Section 7.3 of \cite{Oksendal}. We call the functions  $\sigma $ and $a$ the diffusion coefficients of the process $X_t$ and operator $\mathcal L$ respectively, while $b$ is the drift vector.  Note also that processes with different $\sigma$ may lead to the same PDE as long as $\sigma\cdot\sigma^T$ is the same. Let us  mention    as further  useful  references \normalcolor \cite{EvansSt}, \cite{SalsaBk}, and \cite{Soner}.

\medskip

\noindent {\bf A comment about real-world practice.} If we consider a field of practical application like  quantitative finance, one may ask the question about which of the two known approaches - PDEs versus  martingales and SDEs - is more important in the real practice of derivatives pricing. Here is a partial  answer: since the Black-Scholes equation is a modified form of the Heat Equation, understanding PDEs is very important as a practical tool, see  \cite{WHD}. And the American options add a free boundary problem, a topic that we will find later in the text. But the stochastic approach contains richer information that is useful in understanding aspects of the financial problem, and it can be easier to formulate and compute in some complex models.

\section{\textsf{Fractional diffusion}}
\label{sec.fraclap}

Replacing the Laplacian operator by fractional Laplacians is motivated by the need to represent processes involving anomalous diffusion. In probabilistic terms, it features long-distance interactions instead of the next-neighbour interaction of random walks and the short-distance interactions of their limit, the Brownian motion. The main mathematical models used to describe such processes are the fractional Laplacian operators, since they have special symmetry and invariance properties that makes for a richer theory.  These operators are generators of stable L\'evy processes that include jumps and long-distance interactions. They  reasonably account for observed anomalous diffusion, with  applications in continuum mechanics (elasticity, crystal dislocation, geostrophic flows,...), phase transition phenomena, population dynamics, optimal control, image processing, game theory, finance, and others. See \cite{App, Bert, ContTankov, GiOsh08, MMM, MK2000, VIKH, Woy2001}, see also Section 1.2 of \cite{Vaz14}.

After a very active period of work on problems involving nonlocal operators, there is now well established theory in a number of directions, like semilinear equations and obstacle problems, mainly of stationary type. We are interested here in evolution problems. Instead of the Heat Equation, the basic evolution equation is now
\begin{equation}\label{fract.eq}
\fbox{\large $u_t +(-\Delta)^s u=0$}
\end{equation}
There has been intense work in Stochastic Processes for some decades on this equation, but not in Analysis of PDEs. My interest in the field dates from the year 2007 in Texas in collaboration with Prof. Luis Caffarelli, who was one of the initiators, specially  in problems related to nonlinear diffusion and free boundaries.

It is known that there is well defined semigroup associated with this equation for every $0<s<1$ that solves Cauchy problem \eqref{fract.eq} in the whole space or the typical initial and boundary value problems in a bounded domain, see more below. Though in the limit $s\to 1$ the standard heat equation is recovered, there is a big difference between the local operator $-\Delta$ that appears in the classical heat equation and represents Brownian motion on one side, and the nonlocal family $(-\Delta)^s $, $0<s<1$, on the other side. In the rest of the paper we are going to discuss some of those differences, both for linear and nonlinear evolution equations. We have commented on the origins and applications of the fractional Laplacian and other nonlocal diffusive operators in our previous survey papers \cite{Vaz12Abel} and \cite{Vaz14}.

\subsection{Versions of the fractional Laplacian operator}

Before proceeding with the study of equations, let us  examine the different approaches and defining formulas for the fractional Laplacian operator. We assume that the space variable $x\in\ren$, and the fractional exponent is~$0<s<1$.

\begin{itemize}[leftmargin=*]
\item {\sc Fourier approach}.
First, we may consider the pseudo-differential operator given by the Fourier transform:
\begin{equation}
\widehat{(-\Delta)^s u}(\xi) = |\xi|^{2s}\widehat{u}(\xi)\,.
\end{equation}
This allows to use the very rich theory of Fourier transforms, but is not very convenient for nonlinear analysis which is our final goal. Due to its symbol $|\xi|^{2s}$, the fractional Laplacian can be viewed as a  symmetric differentiation operator of fractional order $2s$. Even when $2s=1$, it is not the standard first derivative, just compare the Fourier symbols.

\item {\sc Hyper-singular integral operator}. The formula reads
\begin{equation}
(-\Delta)^s u(x)= C_{N,s}\int_{\ren}\frac{u(x)- u(y)}{|x - y|^{N+2s}}\,dy\,.
\end{equation}
The kernel is not integrable near $x$ and this motivates the need for the difference in the numerator of the integrand. The integral is understood as principal value.
With this definition, the operator is the inverse of the Riesz integral operator $(-\Delta)^{-s}u,$ which has a more regular kernel $C_1|x - y|^{N-2s}$, though  not integrable at infinity.  The fractional Laplacian operator is also called the {\sl Riesz derivative}.

\item {\sc Numerics and stochastic approach.} Take the random walk for a processes with probability $u^n_j$ at the site $x_j$ at time $t_n$:
\begin{equation}
u^{n+1}_j =\sum_k P_{jk}u^n_k\,,
\end{equation}
where $\{P_{jk}\}$ denotes the transition function which has a {\sl fat  tail} (i.e., a power decay with the distance $|j-k|$), in contrast to the next-neighbour interaction of random walks. In a suitable limit of the space-time grid you get an operator $A$ as the infinitesimal generator of a L\'evy process: if $X_t$ is the isotropic $\alpha$-stable L\'evy process we have
\begin{equation}
Au(x)=\lim_{h\to 0}\frac1{h}\mathbb{E}(u(x)-u(x+X_h)).
\end{equation}
The set of functions for which the limit on the right side exists (for all
$x$) is called the domain of the operator. We arrive at  the fractional Laplacian with an exponent $s\in (0,1)$ that depends on the space decay rate of the interaction $|j-k|^{-(N+2s)}$, $0<s<1$.

\item {\sc The Caffarelli-Silvestre extension.} The $\alpha$-harmonic extension: Find first the solution of the $(N+1)$-dimensional elliptic problem
\begin{equation}
\nabla\cdot(y^{1-\alpha}\nabla U) = 0  \quad (x, y)\in \ren\times\re_+; \quad
U(x, 0) = u(x), \quad  x\in \ren.
\end{equation}
The equation is degenerate elliptic but the weight belongs to the Muckenhoupt $A_2$ class, for which a theory exists \cite{FKS}. We may call $U$ the extended field. Then, putting $\alpha=2s$ we have
\begin{equation}
(-\Delta)^{s} u(x)= -C_{\alpha}\lim_{y\to 0}y^{1-\alpha}
\frac{\partial U}{\partial y}\,.
\end{equation}
When $s= 1/2$, i.e. $\alpha=1$, the extended function $U$ is harmonic (in $N+1$ variables) and the operator is the Dirichlet-to-Neumann map on the base space $x\in \ren$. The general extension was proposed in PDEs by Caffarelli and Silvestre \cite{Caffarelli-Silvestre}, 2007, see also \cite{SilvTh}.  This construction is generalized to other differential operators, like the harmonic oscillator, by  Stinga and Torrea, \cite{StingaTorrea2010}.

\item {\sc Semigroup approach.} It uses the following formula in terms of the heat flow generated by the Laplacian $\Delta$:
\begin{equation}\label{sLapl.semigrp}
\displaystyle(-\Delta)^{s}
f(x)= \frac1{\Gamma(-s)}\int_0^\infty
\left(e^{t\Delta}f(x)-f(x)\right)\frac{dt}{t^{1+s}}.
\end{equation}
\end{itemize}

Classical  references for analysis background on the fractional Laplacian operator in the whole space: the books by N. Landkof \cite{Landkof} (1966), E. Stein \cite{Stein} (1970), and E. Davies \cite{Davies1} (1996). The recent monograph \cite{BucVal16} by Bucur and Valdinoci (2016) introduces fractional operators and more generally nonlocal diffusion, and then goes on to study a number of stationary problems.  Numerical methods to calculate the fractional Laplacian are studied e.g. in \cite{Poz2016}.

\medskip

{\bf Fractional Laplacians on bounded domains}

All the previous versions are equivalent when the operator acts in $\ren$. However, in order to work in a bounded domain $\Omega\subset \ren$ we will have to re-examine all of them. For instance, using the Fourier transform makes no sense. Two main efficient alternatives are studied in probability and PDEs, corresponding to different way in which the information coming from the boundary and the complement of the domain is to be taken into account. They are called the restricted fractional Laplacian (RFL) and the spectral fractional Laplacian (SFL), and they are carefully defined in Section \ref{sec.bdd.domain}. And there are more alternatives that we will also discuss there.


\subsection{Mathematical theory of the Fractional Heat Equation}

 The basic linear problem is to find a solution $u(x,t)$ of
\begin{equation}\label{fractlin.heq}
u_t+(-\Delta)^{s}u=0\,,\qquad 0<s<1\,.
\end{equation}
We will take $x\in\RR^N$, $0<t<\infty$,  with initial data  $u_0(x)$ defined for $x\in\ren$. Normally, $u_0, u\ge 0$, but this is not necessary for the mathematical analysis. We recall that this model represents the linear flow generated by the so-called L\'evy processes in Stochastic PDEs, where the transition from one site $x_j$ of the mesh to another site $x_k$ has a probability that depends on the distance $|x_k-x_j|$ in the form of an inverse power for $j\ne k$, more precisely, \ $c\,|x_k-x_j|^{-N-2s}$. The range is $0<s<1$. The limit from random walk on a discrete grid to the continuous equation can be read e.\,g. in  Valdinoci's \cite{Valdinoc}.

The solution of the linear equation can be obtained in $\ren$ by means of convolution with the fractional heat kernel
\begin{equation}\label{fract.conv}
 u(x,t)=\int u_0(y)P_t(x-y)\,dy,
\end{equation}
 and the probabilists Blumental  and  Getoor proved in the 1960s \cite{Blumenthal-Getoor},   that
\begin{equation}\label{Pkernel}
P_t(x)\asymp  \frac{t}{\big(t^{1/s}+|x|^2\big)^{(N+2s)/2}}\,.
\end{equation}
Here $ a\asymp b$ means that $a/b$ is uniformly bounded above and below by a constant.
Only in the case $s=1/2$ the kernel is known to be explicit, given precisely by the previous formula up to a constant. Note the marked difference with the Gaussian kernel $G_t$ of the heat equation (case $s=1$). The behaviour as $x$ goes to infinity of the function $P_t$ is power-like  (with a so-called {\sl fat tail}) while $G_t$ has exponential spatial decay, see \eqref{gauss.kernel}. This difference is expected in a theory of long-distance interactions. See more on this issue in \cite{Karch09}.

Rather elementary analysis allows then to show that the convolution formula generates a contraction semigroup in all $L^p(\ren)$ spaces, $1\le p \le \infty$, with regularizing formulas of the expected type
$$
\|u(t)\|_\infty\le C(N,s,p) \,\|u_0\|_p\,\,t^{-N/2sp}\,.
$$
When the data and solutions are not assumed to be Lebesgue integrable,  interesting questions appear. Such  questions  have been solved for the classical heat equation, where it is well-known that solutions exist for quite large initial data, more precisely data with square-quadratic growth as $|x|\to \infty$, see Widder \cite{W2}. The idea is that the convolution formula \eqref{fract.conv} still makes sense and can be conveniently manipulated.

Likewise, we may study the fractional heat equation in classes of (maybe) large functions and pose the question: given a  solution of the initial value problem posed in the whole space $\mathbb{R}^N$, is it representable by the convolution formula? The paper \cite{BPSV} by B. Barrios, I. Peral, F. Soria, and E. Valdinoci, shows that the answer is yes if the solutions are suitable strong solutions of the initial value problem posed in the whole space $\mathbb{R}^N$, they are nonnegative, and the growth in $x$ is no more that $u(x,t)\le (1+|x|)^a$ with $a<2s$.

In the recent paper \cite{BSiV2016} by M. Bonforte, Y. Sire, and the author, we look for optimal criteria. We pose the problem of existence, uniqueness and regularity of solutions for the same initial value problem in full generality. The optimal class of initial data turns out to be the class of locally finite Radon measures $\mu$ satisfying the condition
\begin{equation}
\int_{\ren} (1+|x|)^{-(N+2s)}\,d\mu(x)<\infty\,.
\end{equation}
We call this class $\mathcal{M}_s$. We construct  weak solutions for such data, and we prove uniqueness of nonnegative weak solutions with nonnegative measure data. More precisely, we  prove that there is an equivalence between nonnegative measure data in that class and nonnegative weak solutions, which is given in one direction by the representation formula, in the other one by the existence of an initial trace.  So the result closes the problem of the Widder theory for the fractional heat equation posed in $\ren$.
We then review many of the typical properties of the solutions, in particular we prove optimal pointwise estimates and new Harnack inequalities. Asymptotic decay estimates are also found for the optimal class. Here is the general result in that direction.
We want to estimate the behaviour of the constructed solution $u=P^t\ast \mu_0$ for $t>0$ and prove that it is a locally bounded function of $x$ and $t$ with precise estimates. Here is the main result.

\begin{theorem}\label{Thm.upper.gen} Let $u=S_t\mu_0$ the very weak solution with initial measure $\mu_0\in \mathcal{M}^+_s$ and let $\|\mu_0\|_\Phi:= \int_{\ren}\Phi\,d\mu_0$. There exists a constant $C(N,s)$  such that for every $t>0$ and $x\in\ren$
\begin{equation}\label{quant.est}
u(t,x)\le C\,\|\mu_0\|_\Phi \,(t^{-{N/2s}}+t) (1+|x|)^{N+2s}\,.
\end{equation}
\end{theorem}
\noindent Here $\mathcal{M}^+_s$ are the nonnegative measures in the class $\mathcal{M}_s$ and
$\|\mu_0\|_\Phi$ is the associated weighted norm with weight $\Phi(x)=(1+|x|^2)^{-(N+2s)/2}$.
See whole details in \cite{BSiV2016}, Theorem 7.1. The dependence on $t$ cannot be improved. Under radial conditions a better growth estimate in $x$ is obtained. Construction of self-similar solutions with growth in space also follows.

 \smallskip

\subsection{Other nonlocal diffusive  operators}  $\bullet$ Equation  \eqref{fractlin.heq} is the most representative example of a wide class of equations that are used to describe diffusive phenomena with nonlocal, possibly long-range interactions. We can replace the fractional Laplacian by a L\'evy operator $L$ which is the pseudo-differential operator with the symbol $a = a(\xi)$ corresponding to a certain convolution semigroup of measures, \cite{Karch09}. Popular models that are being investigated are integro-differential operators with irregular or rough kernels, as in \cite{SchwSilv}, where the form is
\begin{equation}
 u_t+ b(x,t)\cdot\nabla u-\int_{\ren}\left(u(x+h,t)-u(x,t)\right)K(x,t,h)\,dh=f(x,t).
\end{equation}
See also \cite{AMRT2010, BdP, ChLD, dPQR2016, FK2013, Silv06, DTEJ2017}, among many other references.

\smallskip

\noindent $\bullet$ A different approach is taken by Nystr\"om-Sande \cite{NySande} and Stinga-Torrea \cite{StTor2015}, who define the fractional powers of the whole heat operator and solve
\begin{equation}
(\partial_t-\Delta)^su(t,x)=f(t,x),\quad\hbox{for}~0<s<1.
\end{equation}
In this equation the random jumps are coupled with random waiting times. The authors find the space-time fundamental solution that happens to be explicit, given by
\begin{equation}
K_s(t,x)=\frac{1}{(4\pi t)^{N/2}|\Gamma(-s)|}\cdot\frac{e^{-|x|^2/4t}}{t^{1+s}}=
\frac{1}{|\Gamma(-s)|t^{1+s}}G(x,t),
\end{equation}
for $x\in\ren$,  $t>0$, where $G$ is the Gaussian kernel. The limits $s\to 0, 1$ are singular. Motivations are given, extension methods are introduced, and regularity results proved.

\normalcolor

\section{\textsf{Nonlinear Diffusion}}
\label{sec.nonlindiff}

The linear diffusion  theory has enjoyed much progress, and is now solidly established in theory and applications. However, it was soon observed that many of the equations modeling physical phenomena without excessive simplification are essentially nonlinear, and its more salient characteristics are not reflected by the linear theories that had been developed, notwithstanding the fact that such linear theories had been and continue to be very efficient for a huge number of applications. Unfortunately, the mathematical difficulties of building theories for suitable nonlinear versions of the three classical partial differential equations (Laplace's equation, heat equation and wave equation) made it impossible to make significant progress in the rigorous treatment of these nonlinear problems until  the 20th century was well advanced. This observation also applies to other important nonlinear PDEs or systems of PDEs, like the Navier-Stokes equations and nonlinear Schr\"odinger equations.

\subsection{Importance of Nonlinear PDEs}
The main obstacle to the systematic study of the Nonlinear PDE Theory was the perceived difficulty and the lack of tools. This is reflected in a passage by John Nash (1958). In his seminal paper \cite{Nash58}, he said

 {\sl The open problems in the area of  nonlinear PDE are very relevant to applied mathematics and science as a whole, perhaps more so that the open problems in any other area of mathematics, and the field seems poised for rapid development. It
seems clear, however, that {\sl fresh methods must be employed...}}

and he continues in a more specific way:

{\sl Little is known about the { existence, uniqueness and
smoothness} of solutions of the general equations of flow for a
viscous, compressible, and heat conducting fluid...}

\smallskip

 This is a grand project in pure and applied science and it is still going on. In order to start the work, and following the mathematical style that cares first about foundations, he set about the presumably humble task of proving the regularity of the weak solutions of the PDEs he was going to deal with. More precisely, the problem was to  prove continuity (H\"older regularity) of the weak solutions of elliptic and parabolic equations assuming the coefficients $a_{ij}$ to be uniformly elliptic (positive definite matrices) but only bounded and measurable as functions of $x \in \ren$.  In a rare coincidence of minds, this was done in parallel by J. Nash  \cite{Nash57, Nash58} and the then very young Italian genius E. De Giorgi \cite{DeGiorgi57}\footnote{Strictly speaking, priority goes to the latter, but the methods were different.}. This was a stellar moment in the History of Mathematics, and the ideas turned out to be ``a gold mine'', in L. Nirenberg's words\footnote{Nash and Nirenberg shared the Abel Prize for 2015.}. The results were then taken up and given a new proof by J. Moser, \cite{Moser60}, who went on to establish the Harnack inequality, \cite{Moser64}, a very useful tool in the sequel.

Once the tools were ready to start attacking Nonlinear PDEs in a rigorous way, it was discovered that the resulting mathematics are  quite different from the linear counterparts, they are often difficult and  complex, they turn out to be more realistic than the linearized models in the applications to real-world phenomena, and finally they give rise to a whole set of new phenomena unknown in the linear world. Indeed, in the last decades we have been shown a multiplicity of new qualitative properties and surprising phenomena encapsulated in the nonlinear models supplied by the applied sciences. Some of them are very popular nowadays, like free boundaries, solitons and shock waves. This has kept generations of scientists in a state of surprise and delight. Nonlinear Science rests now on a firm basis and Nonlinear PDEs are a fundamental part of it.


\subsection{Nonlinear heat flows, nonlinear diffusion}

The general formula for the nonlinear diffusion models in divergence form is
\begin{equation}
\fbox{ $ u_t=\sum \partial_i A_{i}(u,
\nabla u) + \sum  {\cal B}(x,u, \nabla u) \,,$}
\end{equation}
where ${\cal A}=(A_i)$ and ${\cal B}$ must satisfy some so-called structure conditions, the main one is again the ellipticity condition on the function $A(x,u,z)$ as a function of the vector variable $z=(z_i)$.
This general form was already posed as a basic research project in the 1960s, cf.  \cite{Serrin64, ArSerr67}. Against the initial expectations,  the mathematical theory  turned out to be  too vast to admit a simple description encompassing the stated generality. There are reference books worth consulting, like those by Ladyzhenskaya et al. \cite{LSU, LU}, Friedman \cite{Fr64}, Lieberman \cite{Lieb96}, Lions-Magenes \cite{LionsMag}, and Smoller \cite{Smo82} are quite useful introductions. But they are only basic references.

Many specific examples, now considered the “classical nonlinear
diffusion models”, have been investigated separately to understand in detail the
qualitative features and to introduce the quantitative techniques, that
happen to be many and from very different origins and types.

My personal experience with nonlinear models of diffusive type lies in two  areas called respectively { `Nonlinear Diffusion with Free Boundaries'} and {`Reaction-Diffusion PDEs'}.

\subsubsection{Pure Nonlinear Diffusion. The Free Boundary Models}

The work on nonlinear parabolic equations in the mathematical research community to which I  belonged focussed attention on the analysis of a number of paradigmatic models involving the occurrence of free boundaries, for which new tools were developed and tested. A rich  theory originated that has nowadays multiple applications.

\medskip

\noindent  $\bullet$  {\sc The Obstacle Problem.}  This is the most famous free boundary problem and there is a huge literature for it, cf. \cite{Caff98, CaffSa05, Fr82, KSt} and their references. It belongs to the class of stationary problems, connected with elliptic equations, hence further away from our  interests in this paper. Let us only say at this point that a free boundary problem is a mathematical problem in which we want to find the solution of a certain equation (normally, a PDE) as well as the domain of definition of the solution, which is also an unkown of the problem. Typically, there exists a fixed `physical' domain $D$ and the solution domain $\Omega$ that we seek is a subset of $D$, well-determined if we know the free boundary $\Gamma=\partial \Omega\cap D$.

\medskip

Parabolic free boundaries may move in time. They appear in the   `four classical sisters' that we will introduce next:

\smallskip

\noindent $\bullet$  {\sc The Stefan Problem} (Lam\'e and Clapeyron, 1833; Stefan 1880)
 The problem typically describes the temperature distribution in a homogeneous medium undergoing a phase change (like ice and water). The heat equation must be solved in both separate media filling together a certain space $D\subset \ren$, and the separation surface is allowed to move with time according to some  transfer law. The mathematical formulation is thus
$$
SE:
\left \{
\begin{array}{ll}
u_t= k_1 \Delta u \quad & \mbox{for } u>0,\nonumber\\
u_t= k_2 \Delta u \quad & \mbox{for } u  <0.
\end{array}\right. \qquad
TC : \left \{
\begin{array}{l}
u=0, \\
{\bf v}= L(k_1\nabla u_1-k_2\nabla u_2).
\end{array}\right.
$$
(SE) means state equations, valid in the separate domains $\Omega_1=\{(x,t): u(x,t)>0\}$ and $\Omega_2=\{(x,t): u(x,t)<0\}$, which are occupied by two immiscible material phases (typically water for $u>0$ and ice for $u<0$). Of course, if the physical domain $D$ is not the whole space then usual boundary conditions have to be given on the fixed boundary $\partial D$. The main mathematical feature is  the existence of a {\sl free boundary} or {\sl  moving boundary}\footnote{also called {\sl interface} in the literature.} $\Gamma\subset \ren\times \re$ that separates ice from water and there $u=0$, see the monographs\cite{Meir, Rub}. This free boundary $\Gamma$ moves in time and has to be calculated along with the PDE solution $u$, so that suitable extra information must be given to determine it: (TC) means transmission condition that applies at the free boundary $\Gamma$, and ${\bf v}$ is the normal advance speed of $\Gamma$. Physically, this formula is due to the existence of latent heat at the phase transition. We not only  want to determine the location of $\Gamma$ but we want to hopefully prove that it is a nice hypersurface in space-time.

Summing up, the combination of analysis of PDEs and variable geometry is what makes this problem difficult. The correct mathematical solution came only via the weak formulation \cite{KamStefan} that allows to eliminate the geometry in a first step and concentrate in finding the so-called weak solution. The free boundary comes later as the zero level set of the weak solution, and finding it needs some regularity theory.

A simpler version is the {\sc One-phase Stefan problem} where ice is assumed to be at zero degrees, roughly $u=0$ in $\Omega_2$. The free boundary is still there but the mathematical theory is much easier, hence better known.

\smallskip

\noindent $\bullet$ {\sc The Hele-Shaw cell problem.} (Hele-Shaw, 1898; Saffman-Taylor, 1958) The problem is posed in a fixed spatial domain $D\subset \ren$, and consists of finding $\Omega(t)\subset D$ and $u(x,t)$ such that
$$
u>0, \ \Delta u=0 \quad \mbox{in} \quad \Omega(t); \quad u=0, \ {\bf v}=L\partial_n u
\quad \mbox{on} \quad \partial_f\Omega(t).
$$
Here the main unknown of the problem  is  the moving domain $\Omega(t)\subset D$, and $\partial_n u$ denotes normal derivative on the free boundary $\Gamma(t)=\partial_f\Omega(t)$, the part of the boundary of the set $\{x\in D: u(x,t)>0\}$ that lies inside  $D$. Additional conditions are to be given on the part of fixed boundary $\partial D$ bounding $\Omega(t)$. Once $\Omega(t)$ is known, solving the Laplace equation for $u$ is standard; notice that it is nontrivial because of the boundary conditions (sometimes there is a forcing term).

Mathematically, this is a simplified version of the previous model where there is only one phase, and besides the time derivative term disappears from the state equation. This increased simplicity comes together with beautiful analytical properties, some of them related to the theory of conformal transformations and complex variables when working in 2D, see \cite{Howison, Rich}. The Hele-Shaw flow appears in fluid mechanics as the limit of the Stokes flow between two parallel flat plates separated by an infinitesimally small gap, and is used to describe various applied problems. The weak formulation is studied in \cite{EllJ}. There are many examples of moving boundaries with interesting dynamics; thus, a peculiar complex variable pattern exhibiting a free boundary with a persistent pointed angle is constructed  in \cite{KLV95} in 2D. In that example, the free boundary does not move until the pointed angle is broken, which happens in finite time. On the other hand, wider angles move immediately and the free boundary is then smooth.

\smallskip

\noindent $\bullet$  {\sc The Porous Medium Equation}. This is an equation in the nonlinear degenerate parabolic category,
$$
u_t=\Delta u^m, \quad m>1.
$$
The equation appears in models for gases in porous media, underground infiltration, high-energy physics, population dynamics and many others. We will devote a whole section to review the free boundary and other nonlinear aspects of this equation, called PME for short, since it has served so much as a paradigm for the mathematics of nonlinear degenerate diffusion, see \cite{ArBk86, VICM06, Vazpme}. Actually, we  see that the free boundary does not appear in the formulation, but it will certainly appear in the theory. It is a {\sl hidden free boundary}.

The equation can be also considered for exponents $m<1$, called {\sl fast diffusion range}, and further generalized into the class of so-called {\sl filtration equations } $u_t=\Delta \Phi(u),$ where $\Phi$ is a monotone increasing real function. This generality also allows to include the Stefan problem that can be written as a filtration equation with very degenerate $\Phi$:
$$
\Phi(u) = (u - 1)_+ \qquad \mbox{for}  \ u \ge  0, \qquad
\Phi(u) = u \qquad  \mbox{for} \ u < 0, \qquad p>1.
$$

\noindent $\bullet$  {\sc The $p$-Laplacian Equation}. This is another model of nonlinear degenerate diffusion
\begin{equation}\label{plap.eq}
u_t=\mbox{\rm div}\, (|\nabla
u|^{p-2}\nabla u)\,.
\end{equation}
Such a model appears in non-Newtonian fluids, turbulent flows in porous media, glaciology and other contexts. The mathematics of this equation  turn out to be closely related to the PME: existence, regularity, free boundaries, and so on, but there are subtle differences. Here $p>2$ is needed for a free boundary to appear, \cite{DiB93bk, VSmoothing}. Recent interest in the limit cases {\sl $p=1$} (total variation flow, used in image analysis), or {\sl $p=\infty$} (appearing in geometry and transport), \cite{Evans}. On the other hand, the equation can be generalized into the class of equations with {\sl gradient-dependent diffusivity} of the general form
\begin{equation*}
u_t - \nabla\cdot (a(|Du|)Du) = 0, \label{gdd.eq}
\end{equation*}
where $a$ is a nonnegative real function with suitable growth assumptions to ensure degenerate parabolicity.
Another extension is the {\sl  doubly non-linear diffusion equation} of the form
$$
u_t = \nabla\cdot(|D(u^m)|^{p-2} D(u^m)).
$$
Here the diffusivity takes the form $a(u,|Du|)=cu^{(p-1)(m-1)}|Du|^{p-2}$. We use the notations $\nabla u=D u$ for the spatial gradient.

\subsubsection{The Reaction Diffusion Models}

This is another important direction taken by Nonlinear Diffusion, in which the nonlinear features originate from a lower-order term with super-linear growth. This may create a mathematical difficulty in the form of {\sl blow-up}, whereby a solution exists for a time interval $0<t<T$ and then some norm of the solution goes to infinity as $t\to T$ (the blow-up time). In other cases  the singular phenomenon is extinction (the solution becomes zero every where), or some other kind of singularity formation.

\noindent $\bullet$  {\sc The Standard Blow-Up model:} It is also called the Fujita model
(Kaplan, 1963; Fujita, 1966)
$$
u_t=\Delta u + u^p\qquad p>1\,.
$$
Main feature: If $p>1$ the norm $\|u(\cdot,t)\|_\infty$ of the
solutions may go to infinity in finite time. This depends on the domain and the initial data. For instance, if the space  domain is $\ren$ and the initial function is constant, then blow-up in finite time always happens. Hint: Integrate the  ODE $u_t=u^p$. However, when the data are distributed in space then diffusion and reaction compete and the result is a priori uncertain. This is how a large literature arose. Thus, if the initial data are bell-shaped (like the Gaussian function), the domain is bounded and boundary conditions are zero Dirichlet, then small data will not blow-up and large data will. For other configurations things depend on the exponent $p$:  there exists a critical exponent $p_F$ called the Fujita exponent, such that all positive solutions blow up if $p\in (1,p_F)$. See \cite{Fuj, GalVazBU, Kaplan63, Lev90}.

A number of beautiful blow-up patterns emerge in such evolutions. Galaktionov and the author have constructed  in \cite{GalVaz.combu} a particular one, called the {\sl peaking solution}, that blows up in finite time $T$ at a single point $x_0$  and then continues for later time as a bounded smooth solution, a clear example of the  curious phenomenon called continuation after blow-up. However, the most common situation in reaction-diffusion systems of this diffusive type is {\sl complete blow-up} at time $T$ with no possible continuation (for instance, the numerical approximation goes to infinity everywhere for $t>T$).
The intricate phenomenon of {\sl bubbling} is studied by M. del Pino in another course of this volume \cite{MdPino}.

As an extension of this elementary reaction-diffusion blow-up model there have been studies for many equations of the general form
$$
u_t={\mathcal A}(u)+ f(u,Du)
$$
where ${\mathcal A}$ is a linear or nonlinear diffusion operator, maybe of porous medium or $p$-Laplacian type. The studies also include systems. Some of them are systems of mixed type, one of the most popular ones is the  chemotaxis system, where blow-up has a very interesting form that is still partially understood, \cite{HV}.

\smallskip

\noindent $\bullet$  {\sc The Fisher-KPP model and traveling waves:}  The problem  goes back to Kolmogorov, Petrovskii and Piskunov, see \cite{KPP}, that present the most simple reaction-diffusion equation concerning the concentration $u$ of a single substance in one spatial dimension,
\begin{equation}\label{KPPdim1}
\partial_t u=D u_{xx} + f(u)\,,
\end{equation}
with an $f$ that is positive between two zero levels $f(0)=f(1)=0$. We assume that $D>0$ is constant. The choice $f(u) = u(1-u)$ yields Fisher's equation \cite{Fisher} that was originally used to describe the spreading of biological populations. The celebrated  result says that the long-time behavior of any solution of \eqref{KPPdim1}, with suitable data $0\le u_0(x)\le 1$ that decay fast at infinity, resembles a traveling  wave with a definite speed that can be explicitly calculated. The KPP traveling wave  pattern is one of the most famous dynamic patterns in diffusive phenomena.

 When considering equation \eqref{KPPdim1} in dimensions $N\geq 1$, the problem becomes
\begin{equation}\label{classicalKPP}
u_t-\Delta u=f(u) \quad \text{in }(0,+\infty)\times \ren,
\end{equation}
This case has been studied by Aronson and Weinberger in
\cite{AronsonWeinberger2,AronsonWeinberger}, where they prove the following result.

\noindent  {\bf Theorem.} {\sl Let $u$ be a solution of \eqref{classicalKPP} with $u_0 \neq 0$ compactly supported in $\ren$ and satisfying $0\leq u_0(x)\leq 1$. Let
$c_*=2\sqrt{f'(0)}$. Then,
\begin{enumerate}
\item if $c>c_*$, then $u(x,t) \rightarrow 0$ uniformly in $\{|x|\geq ct\}$ as $t\rightarrow \infty$.
\item if $c<c_*$, then $u(x,t) \rightarrow 1$ uniformly in $\{|x|\leq ct\}$ as $t\rightarrow \infty$.
\end{enumerate}
}

\noindent In addition, problem \eqref{classicalKPP} admits planar traveling wave solutions connecting $0$ and $1$, that is, solutions of the form
$u(x,t)=\phi(x\cdot e+ct)$ with
$$
-\phi'' + c\phi'=f(\phi) \text{ in } \re, \quad \phi(-\infty)=0, \ \phi(+\infty)=1.
$$
This asymptotic traveling-wave behavior  has been generalized in many interesting ways, in particular in nonlinear diffusion of PME or $p$-Laplacian type, \cite{dPJLVjde1991},
\cite{dePabloSanchez98}, \cite{AuVz16}.   Departing from these results, King and McCabe examined in  \cite{KingMcCabe} a case of fast diffusion,
namely
$$
u_t=\Delta u^{m}  +u(1-u),\quad x \in \mathbb{R}^N, t>0,
$$
where $(N-2)_+/N<m<1$, and  showed  that the problem does not admit traveling wave solutions and the long time behaviour is quite different. We will return to this question when dealing with fractional nonlinear diffusion  in the work \cite{StanVazquezKPP}, in Section \ref{sec.mod2}.

\medskip

\noindent $\bullet$ In the last decades many other models and variants of diffusive systems have been proposed, in particular in the form of systems, like the various {\sl cross-diffusion systems} \cite{Jungel}. Cross-diffusion gives rise to instabilities that attract much attention in  population dynamics, since they allow to predict important features in the study of the spatial distribution of species. The seminal work in this field is due to Alan Turing \cite{Turing52}.
In order to understand the appearance of certain patterns in nature with mathematical regularities like the Fibonacci numbers and the golden ratio, he proposed a model consisting of a system of reaction-diffusion equations.

\medskip

\noindent $\bullet$ Blow-up problems have appeared in related disciplines and some of them have attracted in recent times the attention of researchers for their difficulty and relevance. We present two cases, a  case still requiring more work, a case enjoying big success. The combination of diffusion with nonlinear reaction is in both cases very intricate and leading to the deepest mathematics.

\smallskip

\noindent {\sc The fluid flow models:}  The {\sl Navier-Stokes} or
{\sl Euler} equation systems for incompressible flow. The nonlinearity is quadratic and affects first order terms. Progress is still partial. There is also much work on the related topic of geostrophic flows. We will not enter into more details of such a relevant topic that has a different flavor.

\smallskip

\noindent  {\sc The geometrical models:} The {\sl  Ricci flow} describes the motion of the metric tensor of a Riemannian manifold by means of the Ricci matrix: \
$\partial_t g_{ij}=-2R_{ij}$.  This is a nonlinear reaction-diffusion system, even if this information is not clear in the succint formula. Posed in the form of PDEs by {R. Hamilton,} 1982, it has become a Clay Millenium Problem. Its solution  by {G. Perelman} in 2003 was one of biggest success stories of Mathematics in the 21st century, see \cite{ChowLN2006, MTian2007}. One of the main points in the proof is the study of the modes of blow-up of this system. In order to see that the evolution system of the Ricci flow is a type of nonlinear diffusion, it is convenient to recall the much simpler case of  two-manifolds, since in that case it reduces to a type of fast diffusion called logarithmic diffusion, see below.

Let us finally mention the equations of movement by curvature to the list of geometrical models. Enormous progress has been made in that topic.

\medskip

\noindent {\sc Basic reading for this chapter:}  On Nonlinear Diffusion: \cite{DiB93bk, DaskKenig, Vazpme}. On free boundaries \cite{Diaz85, Fr82}. Moreover, \cite{CrankD, CrankFB}. Fully nonlinear equations are form a vast topic that we have not touched, see \cite{CabCaff}.


\section{\textsf{PME: degenerate diffusion and free boundaries}}
\label{sec.pme}

A very simple model of nonlinear diffusion in divergence form is obtained by means of
the equation
\begin{equation}\label{pme.gen}
\fbox{ $ u_t=\nabla\cdot (D(u)\nabla u)\ $}
\end{equation}
where $D(u)$  is a diffusion coefficient that depends on the `concentration variable' $u$. Strict parabolicity requires that $D(u)>0$, and the condition can be relaxed to degenerate parabolicity if we make sure that $D(u)\ge 0$. Now, if we further assume that $D(u)$ is a power function, we get the simplest model of nonlinear diffusion equation in the form
\begin{equation}\label{pme1}
u_t=\nabla\cdot (c_1|u|^{\gamma}\nabla u)=c_2\Delta (|u|^{m-1}u).
\end{equation}
with $m=1+\gamma$ and $c_1,c_2>0$. Exponent $m$, in principle positive, will play an important role in the model, but the constants $c_i$  are inessential, we may put for instance $c_2=1$, $c_1=m$.
The {\sl  concentration-dependent diffusivity} is then
$$
D(u)=m|u|^{m-1}\,.
$$
In many of the applications $u$ is a density or concentration, hence essentially nonnegative, and then we may write the equation in the simpler form
\begin{equation}\label{pme1p}
u_t=\Delta (u^{m})\,,
\end{equation}
that is usually found in the literature (we have dispensed with useless constants). But there are applications in which $u$ is for instance a height that could take negative values, and then version \eqref{pme1} is needed, since otherwise $D(u)$ would not be positive and the equation would not be  parabolic.

The equation has enjoyed a certain popularity as a mathematical model for degenerate nonlinear diffusion, combining interesting and varied applications with a rich mathematical theory. The theory has many interesting aspects, like functional analysis in the existence and uniqueness theory, and geometry in the study of the free boundaries, as well as deep novelties in the long time asymptotics. Our monograph \cite{Vazpme} gathers a large part of the existing theory up to the time of publication  (2007). We will devote this section to review some of the main topics that affect the theory of fractional porous medium models of later nonlocal sections, and we will also present the very recent sharp results on the regularity and asymptotic behaviour of free boundaries, obtained in collaboration with Kienzler and Koch in \cite{KKV16}.

\subsection{The Porous Medium Equation} As we have already said, the value of exponent $m$ is an important part of the model. Clearly, if $m=1$ we have  $D(u)=1$, and we  recover the classical heat equation, $u_t=\Delta u$,  with its well-known properties, like the maximum principle, the $C^\infty$ regularity of solutions, and the infinite speed of propagation of positive disturbances into the whole space, as well as the asymptotic convergence to a Gaussian profile for  suitable  classes of initial data.

\noindent $\bullet$  The first interesting nonlinear case is $m> 1$ where $D(u)$  degenerates at the level $u=0$. This brings as a consequence the existence of weak solutions that have compact support in the space variable for all times, though that support expands. We refer to that situation as {\sl Slow Diffusion}. As a consequence, free boundaries arise and a whole geometric theory is needed. All this is in sharp contrast with the heat equation.

The differences with the heat equation can be seen by means of an easy calculation for $m=2$. In that case, and under the assumption that $u\ge 0$, the  equation can be re-written  as
 $$
\frac 12 u_t= u\Delta u + |\nabla u|^2\,,
 $$
and we can immediately see that for values $u\gg 0$ the equation looks like a harmless nonlinear perturbation of the heat equation plus a lower order term, while for $u\sim 0$ the first term disappears and the equation looks like (a singular perturbation of) the
eikonal equation
 $$
u_t=|\nabla u|^2\,.
 $$
This last equation is not {\sl parabolic}, but {\sl hyperbolic}, with propagation along characteristics. The PME equation is therefore of mixed type near the critical value $u=0$ where it degenerates, and it has therefore mixed properties.

The calculation may look very particular, for a specific value of $m$. But to the initial surprise of researchers, it extends to all value $m>1$, of the slow diffusion range. The {\sl pressure  transformation} $v=cu^{m-1}$ allows us to get an equivalent equation for $v$:
\begin{equation}\label{press.eq}
v_t=(m-1)v\Delta v + |\nabla v|^2\,,
\end{equation}
where we have used the standard normalization  $c=m/(m-1)$. Indeed, the apparent generality of this transformation goes further and there is a great unity in the theory developed for the PME in the whole range $m>1$, see \cite{Vazpme}.  Indeed, it can be proved that in some weak sense the eikonal equation holds on the free boundary $\{u=0\}$, and this implies that the support of the solution spreads with time, another property that can be rigorously proved.

\medskip

$\bullet$ The pressure transformation is even more general, and can be applied to the filtration equation $u_t=\Delta \Phi(u).$ If we put \ $v=\int_1^u (\Phi'(s)/s)\, ds$,  then we can get the pressure equation
\begin{equation}\label{press.eq2}
v_t=\sigma(v)\Delta v + |\nabla v|^2\,,
\end{equation}
where the function $\sigma(v)=\Phi'(u)\ge 0$, cf. \cite{BrnVaz}.

\medskip

$\bullet$ These pressure considerations apply under the assumption that $u\ge 0$, which is physically natural for most applications. It must be pointed out the existence and uniqueness theory has been done for signed solutions, according to the generality that is suitable in Functional Analysis. However, many of the estimates on which the qualitative theory is based do not apply for general signed solutions, and we will forsake them and assume $u\ge 0$ in the rest of the section unless mention to the contrary.

\medskip

$\bullet$ When $m<1$ the equation becomes  singular at $u=0$ in the sense that $D(u)\to\infty$. This range is called {\sl Fast Diffusion}. We will return to that case in the next Section \ref{subsec.fde} since its properties show a remarkable difference with the PME range $m>1$.


\subsection{Applied motivation. Fixing some physical concepts}

This application is maybe the best known and has played a role in developing the theory for the PME, a clear example of the influence of physics on the mathematics. According to Leibenzon (1930) and Muskat (1933), the flow of gas in a porous medium (they were thinking of the petroleum industry) obeys the laws
$$
\left\{\begin{array}{c}
\rho_t + \mbox{\rm div}\, (\rho {\bf V})=0,\\[6pt]
{\bf V}=-\frac{k}{\mu}\nabla p, \quad p=p(\rho),
\end{array}\right.
$$
where $\rho$ is density, $p$ is the averaged pressure and $\bf V$ is the seepage velocity. The first line is the usual continuity equation from fluid mechanics, and the second line left is the {\sl Darcy law} for flows in porous media (Darcy, 1856). Therefore, these porous media flows are potential flows due to averaging of Navier-Stokes on the pore scales. We need a precise closure relation which is given by a gas law of the form
 $p=p_o \,\rho^{\gamma}$, with value of the exponent $\gamma=1$
(isothermal gas) or  $\gamma>1$ (adiabatic gas flow), see details in \cite{Vazpme}. Hence, we get
$$
\rho_t= \mbox{\rm div}\, (\frac {k}{\mu} \rho \nabla p)= \mbox{\rm
div}\, (\frac {k}{\mu} \rho \nabla (p_o \rho^{\gamma})) =c\Delta
\rho^{\gamma+1}.
$$
In order to get the PME we put $u=\rho$, $m=1+\gamma$ (which happens to be equal or larger than 2) and we eliminate useless constants. We point out that the pressure is then
$$
p=p_o u^{m-1},
$$
just the variable that we called $v$ in formula \eqref{press.eq}. No wonder that this equation is important. As for the local flow velocity we have ${\bf V}=-c\nabla v$ in our mathematical notation.

- There are many other applications, as described in the book \cite{Vazpme}: underground water infiltration (Boussinesq, 1903) with  $m=2$, plasma radiation with $m\ge 4$, (Zeldovich-Raizer, around 1950), spreading of populations (self-avoiding diffusion)
{\sl $m\sim 2$}  (Gurtin-McCamy, 1977), thin films under gravity with no surface tension
{\sl $m=4$}, and so on.

\subsection{Generalities. Planning of the theory}

The way the nonlinear theory of the PME has developed is quite different from the way the linear heat equation is studied. Indeed, in the early years there were attempts to construct a perturbation theory putting $m=1+\ve$ in \eqref{pme1p} and then perturbing the linear model, but the singular perturbation analysis was not successful. Fortunately, around  1958 when the theory started the serious development in Moscow \cite{OKC58}, the tools of nonlinear functional analysis were ready, and in particular the concept of weak solution and the role of a priori estimates.

Here are the main topics of mathematical analysis (1958-2016):

 - The precise meaning of solution. Since it was realized that classical solutions do not exist if there are free boundaries.

\medskip

 - The nonlinear approach: estimates; functional spaces.

\medskip

 - Existence of suitable solutions (like weak solutions).  Uniqueness. Further in the theory, variant of the equation showed cases of non-existence or non-uniqueness.

\medskip

- Regularity of solutions: Are weak solutions indeed continuous functions? are they $C^k$ for some  $k$? which is the optimal $k$?

\medskip

- Existence, regularity and movement of interfaces: are they
 $C^k$ for some $k$?

\medskip

- Asymptotic behaviour: is there something comparable to the Gaussian profile as a universal attractor? This is a question of {\sl emerging  patterns}. If there is convergence to a pattern we want to know that rate of convergence. We also want to know how universal that convergence is, in other words the basin of attraction of the asymptotic pattern.

\medskip

- Comparison with other approaches like the probabilistic approach. Interesting new tools appear, like Wasserstein metrics and estimates.

\medskip

The beauty of this plan is that it can be used  \sl mutatis mutandis \rm  on a huge number of related models: fast diffusion models, inhomogeneous media, anisotropic
media, $p$-Laplacian models, applications to geometry or image processing; equations involving effects, like the chemotaxis models,...


\subsection{Fundamental solutions. The Barenblatt profiles }

  These profiles are the alternative to  the Gaussian profiles of the linear diffusion case.
They are source solutions. {\sl  Source} means that $u(x,t)\to M\,\delta(x)$ as  $t\to 0$.
 Explicit formulas exist for them (1950):
\begin{equation}\label{Barensol}
{\bf B}(x,t;M)= t^{-\alpha} {\bf F}(x/t^{\beta}), \quad {\bf
F}(\xi)=\left(C - k \xi^{2}\right)_+^{1/(m-1)}
\end{equation}

\hspace{-1cm}
\includegraphics[width=6.5 cm]{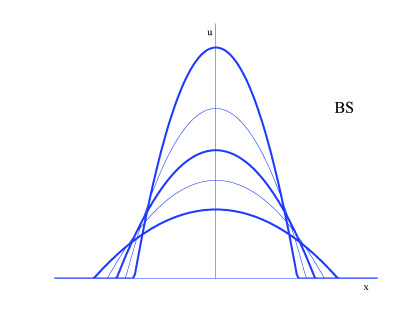}
\vspace{-4cm}
$$\hspace{6cm} \begin{array}{l}
               \alpha= \frac{n}{2+N(m-1)}\vspace{3mm}\\
              \beta=\frac{1}{2+N(m-1)}<1/2\vspace{3mm}\\
                 \mbox{\rm  Maximum height $u=Ct^{-\alpha}$}\vspace{3mm}\\
                 \mbox{\rm  Free boundary is located at distance
                 $|x|=ct^{\beta}$}
                \end{array}
              $$

\

\noindent where $C>0$ is a free constant and $k=k(m,N)$. Since Fourier analysis is not a way for find them, new ideas are needed. We observe that the solution obeys a scaling symmetry, it is self-similar. In other words, it is invariant under suitable scaling in space and time. This fact is the key to finding the expression, see \cite{Vazpme}, page 63.
An important property for the applications is that  \eqref{Barensol}  breaks with the Brownian space-time law: $|x|=ct^{1/2}$,  so that it can be classified as {\sl anomalous diffusion}.

If you look for the mathematical properties, we find a surprise with regularity.
Put $m=2$ for simplicity. ${\bf B}(x,t;M)$ does not satisfy the equation in a classical sense since $u$ is not even $C^1$ continuous in space or time.
The validity of this physical solution was a hot problem when it was discovered around 1950.


\subsection{Concepts of solution}

Hence, there is a problem with the concept of solution that will satisfy the mathematical requirements (existence and uniqueness for a reasonable class of data, plus stability estimates) as well as the physical requirements (to reflect the behaviour that is expected from the evidence obtained in the applications). This problem did not exist for the main example of diffusion, the Heat Equation, since classical solutions could be found.

Many concepts of generalized solution have been used in developing the mathematical theory of the PME, and also in many related equations, not only in the parabolic theory:

$\bullet$ {\sl Classical solution.} This is the most desirable option, and indeed it happens for non-degenerate situations, $u>0$. But it cannot be expected if the Barenblatt solutions are to be included.

$\bullet$ {\sl Limit solution.} This is the practical or computational remedy. To replace the equation by approximated problems with good physical or computational properties and then to pass to the limit. The catch is that the approximation may not converge, or we could be unable to prove it;  even if the approximations do converge, the limit may depend  on the approximation. Spurious solutions may appear when the approximation is not efficient, a quality difficult to tell a priori. On the positive side, limit solutions have been successfully used in the diffusive literature with the names of minimal solutions, maximal solutions, SOLAs (solutions obtained as limits of approximations), proper solutions, ...

$\bullet$ {\sl Weak solution.} This was a very good solution to the problem of building a theory for the PME. The idea is to test the equation against a full set of smooth functions and to eliminate all or most of the derivatives prescribed by the equation on the unknown function. It was first implemented on the PME by O. Oleinik and collaborators \cite{OKC58} (1958). The simplest weak version reads
$$
\int\int (u\,\eta_t - \nabla u^m\cdot\nabla \eta)\,dxdt+\int
u_0(x)\,\eta(x,0)\,dx=0\,,
$$
while there is a second version, the {\sl very weak solution},
$$
\int\int (u\,\eta_t + u^m\,\Delta \eta)\,dxdt+\int
u_0(x)\,\eta(x,0)\,dx=0.
$$
This version is more relaxed than the first. In both cases functional spaces have to be chosen for the solutions to belong to so that the integrals in the formulation make sense and existence and uniqueness can be proved.

Once existence and uniqueness of a weak solution was proved for suitable initial data; that it was verified that all classical solutions are weak; and that the Barenblatt solutions are indeed solutions for $t>0$ (and take the initial Dirac delta in a suitable sense) the theory of the PME could be conveniently framed as  a theory of weak solutions. Even more, it was proved that suitable numerical approximations converge to the weak solution.

$\bullet$  {\sl Better regularity. Strong solution.} The previous paragraph solves the problem of the correct setting in principle. But researchers want to have solutions that have good properties so that we can do calculus with them. Fortunately, weak solutions of the PME are better than weak, they are strong. In this context it means that all weak derivatives entering the original equation are $L^p$ functions for some $p$.

$\bullet$ The search for an abstract method to solve a large number of evolution problems of diffusive type has led to a functional approach called {\sl mild or semigroup solution}, that we discuss below.

$\bullet$  Solutions of more complicated diffusion-convection equations have motivated new concepts that can be translated to the PME:

- {\sl Viscosity solution}. {Two different ideas: (1)  add artificial viscosity and pass to the limit; (2) viscosity concept of Crandall-Evans-Lions  (1984);  adapted to PME
by Caffarelli-V\'azquez (1999).}

-  {\sl Entropy solution} (Kruzhkov, 1968).
Invented for conservation laws;  it identifies unique physical
solution from {\sl spurious} weak solutions. It is useful for
general models with degenerate diffusion plus convection.

-  {\sl Renormalized solutions} (by Di Perna - P. L. Lions), {\sl BV solutions } (by Volpert-Hudjaev),  {\sl Kinetic solutions } (by Perthame,...).

\subsection{Semigroup approach. Mild solution}

Functional Analysis is a power tool for the expert in PDEs, and when used wisely it produces amazing result. Thus, when faced with the task of solving evolution equations of the type
\begin{equation}
u_t + A u=0\,,
\end{equation}
where $A$ is a certain operator between function spaces, we may think about discretizing the evolution in time by using a mesh $t_0=0<t_1,... t_K=T$ and posing the implicit problems
$$
\frac{u(t_k)-u(t_{k-1})}{h_k} + A(u(t_k))=0,\quad h_k=t_k-t_{k-1}.
$$
In other words, we want to find a discrete approximate solution  $u=\{u_k\}_k$ such that
\begin{equation}\label{itd.k}
h A(u_k)+ u_k= u_{k-1}\,,
\end{equation}
where we have used equal time spacing $h_k=h>0$ for simplicity. Of course, the approximate solution depends on the time step $h$, so that we should write $u^{(h)}=\{u_k^{(h)}\}_k$.
This step is called {\sl Implicit Time Discretization,} ITD. We start the iteration by assigning the initial value $u_0^{(h)}=u_{0h}$, where $ u_{0h} $ is the given initial data or an approximation thereof.

\medskip

\noindent {\bf Parabolic to Elliptic}. The success of ITD depends on solving the iterated equations \eqref{itd.k} in an iterative way. In fact, the iteration has always the same format
\begin{equation}\label{itd.st}
h A(u)+ u= f\,,
\end{equation}
since $f=u_{k-1}$ is the value calculated in the previous step. When this is used for the filtration equation $u_t -\Delta \Phi( u)=0,$ we get the stationary equation
\begin{equation}\label{itd.st2}
-h \Delta \Phi(u)+ u= f\,,
\end{equation}
and the question reduces to solve for $u$ if $f$ is known. An easy change of variables
$v=\Phi(u)$, $u=\beta(v)$ leads to
\begin{equation}\label{itd.st.v}
-h \Delta v + \beta(v)=f.
\end{equation}
This is the semilinear elliptic problem that we must solve. We have reduced the theory of a (possibly nonlinear) parabolic problem to an elliptic problem with a specific form.

\medskip

\noindent {\bf Accretive operators. Semigroup generation.}
The rest of the story depends on the theory of accretive operators. If $A$ is an $m$-accretive map in a Banach space $X$ with densely defined domain, then the famous
Crandall-Liggett Theorem  \cite{CL71} ensures not only existence of the solution of the iterated problems, but  also that as $h\to 0$ the discrete solutions $u^{(h)}$ converge to a function $u(t)\in C([0,\infty):X)$ that solves the evolution problem in a sense called {\sl mild sense}. The solution is often termed the {\sl semigroup solution}. Moreover, the set of solutions forms a semigroup of contractions in $X$.

But is this mild solution a solution in some more usual sense? In the case of the PME it is proved that the operator given by $A(u)=-\Delta (u^m)$ is $m$-accretive in the space $X=L^1(\ren)$ when properly defined, \cite{BBC75}, and also that the mild solution is a weak solution. We have explained the method in some detail in Chapter 10 of \cite{Vazpme}.

\subsection{Regularity results for nonnegative solutions}

The next step in the theory is proving that under mild conditions on the data weak solutions of the PME are indeed continuous, and the free boundary is quite often a regular hypersurface in space-time, or in space for every fixed time. We recall that we are working with nonnegative solutions. We we also dealing with the Cauchy problem in the whole space to save effort and concentrate on the basics, but many results hold for locally defined solutions.

$\bullet $ The regularity theory for solutions relies on the existence of a rather miraculous a priori estimate, called the Aronson-B\'enilan estimate \cite{AB79}, that reads:
$$
\Delta v\ge -C/t,
$$
where $v=c u^{m-1}$ is the pressure and $ C=(n(m-1)+2)^{-1}$. Nonnegative solutions with data in any $L^p$ space are then proved to be bounded for positive times. A major step was then done by Caffarelli and Friedman (1982) when they proved  $C^\alpha$ regularity: there is an $\alpha\in (0,1)$ such that a bounded solution defined in a cube is $C^\alpha$ continuous. This holds in all space dimensions.

$\bullet $ What happens to the free boundary? It was soon proved that free boundaries may be stationary for a while but eventually they must move to fill the whole domain as time passes. The movement is expansive, the positivity set keeps expanding in time and never recedes.  Caffarelli and Friedman proved subsequently that if there is an interface $\Gamma$, it is also a $C^\alpha$ continuous set in space time (properly defined).

$\bullet $ How far can you go? The situation is understood in 1D. On the one hand, free boundaries can be  stationary for a time (metastable) if the initial profile is quadratic near $\partial\Omega$: $v_0(x)=O(d^2)$, $d$ being distance to the zero set. This time with lack of movement is called a {\sl waiting time.} It was precisely characterized by the author in 1983;  it is visually interesting in the experiments with thin films spreading on a table. In paper \cite{ACV85} we proved that metastable interfaces in 1D may start to move abruptly after the waiting time. This was called a {\sl corner point}. It implies that the conjecture of $C^1$ regularity for free boundaries in 1D was false. But in 1D the problems with regularity stop here: 1D free boundaries are strictly moving and $C^\infty$ smooth after the possible corner points. See \cite{Vazpme} for full details.

\subsection{Regularity of free boundaries in several dimensions}

The situation is more difficult for free boundary behaviour in several space dimensions, and the investigation is still going on.
\begin{center}
\includegraphics[width=3.cm]{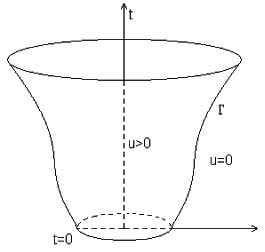}\\
{A regular free boundary in N-D}
\end{center}

$\bullet$ Caffarelli, Wolanski, and the author proved in 1987 that  if $u_0$ has compact support, then after some time $T>0$ the interface is $C^{1,\alpha}$, and the pressure is also $C^{1,\alpha}$ in a lateral sense \cite{CVW87}, \cite{CW90}. Note the lateral regularity is the only option since the Barenblatt solutions are an example of solutions that exhibit a smooth profile that is broken at the FB. The general idea, taken from 1D, is that  when the FB moves, the adjoining profile is always a broken profile, since the support of the solution moves forward only if the gradient of the pressure is nonzero (Darcy's law).

$\bullet$ In his excellent doctoral thesis (1997), Koch  proved that if $u_0$ is compactly supported and transversal then the free boundary is $C^\infty$ after some finite time and the pressure is ``laterally" $C^\infty$. This solved the problem of optimal regularity in many cases, though not all.

$\bullet$  The free boundary for a solution with a hole in 2D, 3D is the physical situation in which optimal regularity can be tested. Indeed, as the flow proceeds and the hole shrinks,  it is observed that the part of the motion of the free boundary surrounding the hole accelerates, so that at the point and time where the hole disappears (this phenomenon is called focusing), the advance speed becomes infinite. The applied setup is a viscous fluid on a table occupying an annulus of radii $r_1$ and $r_2$. As time passes $r_2(t)$ grows out while $r_1(t)$ goes to the origin. In a finite time $ T$ the hole disappears. The flow can be regular for $t<T$ but it was suspected from the numerical evidence that it was not at  the focusing time $t=T$.

To prove this fact, a self-similar solution was constructed displaying the focusing behaviour. It has the form
$$
U_f(x,t)=(T-t)^{\alpha}F(x/(T-t)^{\beta}).
$$
with $(m-1)\alpha = 2\beta+1$. The profile is such that $F(\xi)\sim |\xi|^{\gamma}$ near $\xi=0$, with $\gamma=\alpha/\beta$. There is one free parameter, let us say $\beta$, that is not known a priori. We find a very interesting mathematical novelty,  an {\sl anomalous exponent}, or {\sl similarity of the second kind} in the terminology popularized by Barenblatt \cite{Barbk96}. The problem was solved by ODE analysis in 1993 by Aronson and Graveleau \cite{ArG93}, and then further investigated by Aronson and collaborators, like  \cite{AnAr95, ArV95}. It is proved that $\gamma<1$ so that the  speed ${\bf v}\sim \nabla u^{m-1}$  blows-up at the focusing time $t=T$ at $x=0$. Moreover, the limit profile $U_f(x,T)$ is not Lipschitz continuous at $x=0$, it is only $C^\gamma$ continuous. This is a counterexample to the hypotheses of higher regularity of heat equation and similar diffusive flows: a degenerate equation like PME has limited regularity for nonnegative solutions with moving free boundaries. Such a phenomenon is known to happen in the other typical evolution free boundary problems. Stefan, Hele-Shaw and $p$-Laplacian equation. For the latter see
\cite{AGV98}.

Summing up, higher regularity for PME flows has an obstacle. We may hope to prove higher regularity if we avoid it, like the situation of compactly supported solutions for large times.

\subsection{Recent results on regularity and asymptotics}

The question remained for many years to know if we can prove regularity of the solutions and their free boundaries under some certain geometrical condition on the solution or the data less stringent than the  conditions of compact support and initial transversality of  papers \cite{CVW87}, \cite{CW90}, \cite{Koch}.

\noindent $\bullet$ Much progress was done recently in paper with H. Koch and C. Kienzler \cite{KKV16}, preceded by \cite{Kienz13}.
Here is the main theorem proved in \cite{KKV16} about regularity of  solutions that are locally small perturbations from a flat profile.

\begin{theorem}\label{th-kvv1} There exists $ \delta_0>0$ such that the following holds: \\
  If $ u $ is a nonnegative $ \delta $-flat solution of the PME at
  $ (0,0) $ on scale $ 1 $ with $ \delta $-approximate direction $ e_n
  $ and $ \delta $-approximate speed $ 1 $, and $ \delta \leq
  \delta_0$, then for all derivatives we have uniform estimates
\begin{equation}
 | {\partial_t}^k {\partial_x}^\alpha \nabla_x (u^{m-1} -(x_n+ t)) | \le C \delta
\end{equation}
at all points $(t,x)\in ([-1/2,0] \times \overline{B_{1/2}(0)})\cap
\mathcal P(u) $ with $C=C(N,m,k,\alpha)>0$. In particular,
$\rho^{m-1}$ is smooth up to the boundary of the support in
$(-\frac1{2},0]\times B_{1/2}$, and
\begin{equation}
|\nabla_x  u^{m-1}-e_n |, \quad |\partial_t u^{m-1} - 1|  \le C \delta\,.
\end{equation}
Moreover, the level sets for positive values of $u$ and the free boundary are uniformly smooth hypersurfaces inside $(-\frac1{2} ,0]\times B_{\frac12}(0)$.
\end{theorem}

The technical assumption is being  $ \delta $-flat, which means being very close to a flat travelling wave (the special solution that serves as model) in a certain space-time neighbourhood. See Definition 1 in the paper. The size of this neighbourhood is taken to be unit, but this is not a restriction by the scale invariance of the equation. The very detailed form of the estimates allows us then to derive very strong results for large times, that we will explain in next subsection.

\noindent $\bullet$  Theorem \ref{th-kvv1} implies the eventual $C^\infty$-regularity
result for global solutions that we were looking for.  The following result is Theorem 2 of \cite{KKV16}. We use the notation $R_B(t)=c_1(N,m)\,M^{(m-1)\lambda}t^{\lambda}$, with $\lambda=1/(N(m-1)+2)$, for the Barenblatt radius for the solution with mass $M$ located at the origin.

\begin{theorem}\label{thm2} Let $u \ge 0$ be a solution of the PME
  posed for all $x\in \ren$, $N\ge 1$, and $t>0$, and let the initial
  data $u_0$ be nonnegative, bounded and compactly supported with
  mass $M= \int u_0 dx>0$.  Then, there exists a time $T_r$
  depending on $u_0$ such that for all $t>T_r$ we have:

(i) {\rm Regularity.} The pressure of the solution $u^{m-1}$ is a $C^\infty$
  function inside the support and is also smooth up to the free
  boundary, with $\nabla u^{m-1} \ne 0 $ at the free
  boundary. Moreover, the free boundary function $t=h(x)$ is
  $C^\infty$ in the complement of the ball of radius $R(T_r)$.

(ii) {\rm Asymptotic approximation.} There exists $c>0$ such that
\begin{equation} \label{precise}
  \begin{split}
    t^{-N\lambda} \Big(a^2 M^{2(m-1)\lambda}- c t^{-2\lambda} -
    \frac{\lambda |x-x_0|^2}{2 t^{2\lambda}} \Big)^{\frac{1}{m-1}}_+
    \le & u(t,x)\\ & \hspace{-4cm} \le t^{-N\lambda} \Big(a^2
    M^{2(m-1)\lambda}+ c t^{-2\lambda} - \frac{\lambda |x-x_0|^2}{2
      t^{2\lambda}} \Big)^{\frac{1}{m-1}}_+
\end{split} \end{equation}
where   $x_0=M^{-1}\int x\rho(x)dx$  is the conserved center of mass, and $a$ is a certain constant. Moreover,
\begin{equation} \label{precisesupport}
B_{R_B(t)-c t^{-\lambda}}(x_0) \subset {\rm supp}(u(\cdot,t)) \subset B_{R_B(t)+ct^{-\lambda}}(x_0)
\end{equation}
\end{theorem}

In this way we are able to solve the problem posed in 1987, and improved by Koch in 1997.
We use delicate flatness conditions, scalings, heat semigroups  and harmonic analysis. We have eliminated the non-degeneracy condition on the initial data.
The estimates are uniform. The result cannot be improved in a number of directions. Besides, some more information is available: if the initial function is supported in the ball
  $B_R(0)$, then we can write the upper estimate of the regularization
  time as
\begin{equation}\label{lowerboundontime}
    T_r=T(N,m) M^{1-m}  R^{\frac1\lambda}.
\end{equation}
By scaling and space displacement we can reduce the proof to the case
$M=1$ and $x_0=0$. The fine asymptotic analysis uses also the results of Seis \cite{Seis1}.

\medskip

\noindent $\bullet$ { \bf Nonlinear Central Limit Theorem revisited}

The last part of Theorem \ref{thm2} refers to the way a general solution with compact support approaches the Barenblatt solution having the same mass. This kind of result is what we have called the PME version of the Central Limit Theorem.

It was proved in due time that the standard porous medium flow has an asymptotic stabilization property that parallels the stabilization to the Gaussian profile embodied in the classical Central Limit Theorem if we take as domain $\mathbb{R}^N$ and data  $u_0(x)\in L^1(\mathbb{R}^n)$. The convergence result is
\begin{equation}\label{af1.pme}
\|u(t)-B(t)\|_1\to 0
\end{equation}
as $t\to \infty$, as well as
\begin{equation}\label{af3.pme}
t^{N\lambda}\,|u(x,t)-B(x,t)|\to 0,
\end{equation}
uniformly in $x\in \mathbb{R}^N$. Here, $B(x,t;M)$ be the Barenblatt with the asymptotic mass $M$. Note that the factor $t^{N\lambda}$ is just the normalization needed to work with relative errors since $B(x,t)$ decays like $O(t^{-N\lambda})$. Proofs are due to Kamin and Friedman \cite{FrKa80} for compactly supported solutions, and the author (2001) in full generality. The result is reported with full detail in \cite{Vazpme} and explained in \cite{VICM06, Vaz14}.

An improvement of the result to indicate a definite rate of convergence is due to Carrillo and Toscani  (2000). It works for solutions with a finite second moment, $\int u_0(x)\,|x|^2\,dx$, \cite{CaTo00} and uses the powerful machinery of entropy methods, that become subsequently very popular in studies of nonlinear diffusion.

The result that we obtain above points out to a finer  error rate for compactly supported solutions, that can be written as
\begin{equation}
t^{N(m-1)\lambda}\,|u^{m-1}(x,t)-B^{m-1}(x-x_0,t;M)|=O(t^{-2\lambda})\,.
\end{equation}
Seis' analysis and our paper show optimality of this rate. Note that $2\lambda <1$.

\medskip

\noindent {\bf Other problems.} There are numerous studies of the PME in other settings, like bounded domains with Dirichlet or Neumann conditions, PME with forcing term:
$$
u_t=\Delta u^m+ f,
$$
PME with variable coefficients or weights, generalized filtration equation, PME with convection and/or reaction, ...

Many of the above results have counterparts for the $p$-Laplacian flow. Thus, stabilization to the $p$-Laplacian version of the Barenblatt solution is proved by Kamin and the author in \cite{KVplap}. There are many studies but no comparable fine analysis of the FB has been done.

\smallskip

\noindent {\bf Further reading for this chapter:} On the PME: \cite{Vazpme, VICM06}. On asymptotic behaviour: \cite{Vaz03} and \cite{Vaz04bdd}. About estimates and scaling: \cite{VSmoothing}. For entropy methods \cite{AMTU, CaTo00, CJMU2001}.

\section{\textsf{The Fast Diffusion Equation}}
\label{subsec.fde}

We will consider now that range $m<1$ for the model equation \eqref{pme1}. The equation becomes  singular at the level $u=0$ in the sense that $D(u)\to\infty$. This range is called {\sl Fast Diffusion Equation}, FDE, and we also talk about {\sl singular diffusion}. The new range  was first motivated by a number of applications to diffusive processes with fast propagation:  plasma Physics (Okuda-Dawson law) \cite{BHoll78}, material diffusion (dopants in silicon) \cite{King88}, geometrical flows (Ricci flow on surfaces and the Yamabe flow), diffusive limit of kinetic equations, information theory, and others, see \cite{VSmoothing}.
Once the mathematics started, it was seen that the FDE offers many interesting mathematical challenges, and some unexpected connections with other disciplines like Calculus of Variations.

The common denominations {\sl slow} and {\sl fast} for the parameter ranges $m>1$ and $m<1$ in \eqref{pme1} refer to what happens for $u\approx 0$. But when large values of $u$ are  involved, the names are confusing since the situation is reversed:

- $D(u)\to \infty $ as $u\to \infty$ if $m>1$ (``slow case'')

- $D(u)\to 0 $ as $u\to \infty$ if $m<1$ (``fast case'')

Indeed, power functions are tricky. The pressure transformation can also be used for the FDE range, but then the factor $m-1$ changes sign and the equation that we obtain is different because of the sign changes. Putting $v=cu^{m-1}$ with $c=m/(1-m)$ we get
a new pressure equation of the form
\begin{equation}
v_t=(1-m)v\Delta v - |\nabla v|^2\,,
\end{equation}
so that the eikonal term is now an absorption term. Note that  $u\to 0$ implies $v\to\infty$ in the FDE.

 \subsection{Barenblatt solutions in the good range}
We have well-known explicit formulas for source-type self-similar solutions
called  Barenblatt profiles, valid for  with exponents $m$ less than 1, but only if
$1>m> m_*=(N-2)/N$ if $N\ge 3$:
$$
{\bf B}(x,t;M)= t^{-\alpha} {\bf F}(x/t^{\beta}), \quad {\bf
F}(\xi)=(C + k \xi^{2})^{-1/(1-m)}.
$$

\begin{center}
\includegraphics[width=5.0 cm]{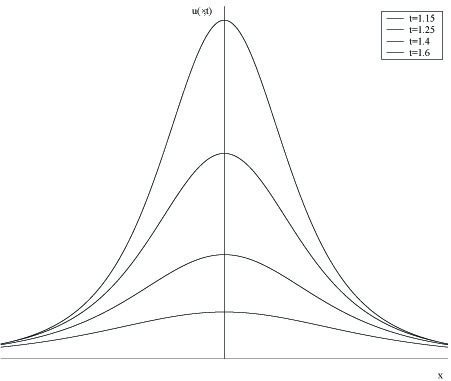}
\end{center}

\noindent The decay rate and spreading rate   exponents are
$$
\alpha= \frac{N}{2-N(1-m)}, \qquad \beta=\frac{1}{2-N(1-m)}>1/2\,.
$$
Both exponents $\alpha, \beta\to \infty$ as $m$ goes down to $m_*$. So the question is what happens for $m<m_*$? It is a long and complicated story, see a brief account further below.

The decay of the Barenblatt FDE profile for fixed time is ${\bf B}=O(|x|^{-2/(1-m)})$, a power-like decay that we have termed a {\sl fat tail} in terms of probability distributions. The exponent ranges from $N$ to $\infty$ in the range $m_*<m<1$.  Note that for $m_*< m <(N-1)/(N+1)$ the distribution ${\bf B}(\cdot,t)$ does not even have a first moment.

The exponent range $m_c<m<1<$ of the FDE where the Barenblatt solutions exist is called the ``good fast diffusion range'', since it has quite nice properties; though different from the linear heat equation, they nevertheless quite satisfying from many points of view, in particular from the point of view of existence, functional analysis, regularity and asymptotic behaviour. Thus, existence of a classical $C^\infty$ smooth solution is guaranteed for every nonnegative, locally bounded Radon measure as initial data (no growth conditions like the HE or the PME), the solution is unique, it is also positive everywhere even if the data are not (they must be nontrivial), decay in time of the solutions depends in a predictable way from suitable norms of the data. Even unbounded Borel measures can be taken, see  \cite{ChassVaz}.

 In particular, when $u_0\in L^1(\ren)$, a semigroup of contractions if generated, conservation of mass holds, and the solution converges for large time to the Barenblatt solution given above, and  with the same expression of the Central Limit Theorem that we saw for the PME. The impressive Aronson-B\'enilan estimate is now a two-sided universal estimate: $- C_1u/t \le u_t(x,t) \le C_2 u/t$, which implies  better and easier estimates for the rest of the theory. Of course, the absence of free boundaries makes it lose part of its power of attraction. We will not enter into the proofs of these results, than can be found in the literature, \cite{DaskKenig, DiB93bk}, \cite{VSmoothing},...

\subsection{Comparison of anomalous diffusions}

The type of diffusion described by the Barenblatt solutions is  called  {\sl anomalous diffusion} since it breaks the Brownian spread rate $|x|\sim t^{1/2}$ and space decay with an exponential rate. We have already seen that anomalous behaviour in the linear setting, as the fundamental solution $P_t(x)$ of the Fractional Heat Equation (FHE) $u_t+(-\Delta)^s u=0$, cf. \eqref{Pkernel}.  In that case the spreading rate is $|x|\sim t^{1/2s}$. This leads to  a formal (and partial) equivalence between anomalous diffusion of FHE and FDE types based on the spread strength. It reads \ $2s\sim 2-N(1-m)$, hence
$$
\ N(1-m)\sim 2(1-s),
$$
which agrees for the classical heat equation: $m=1$, $s=1$. For the best known case $s=1/2$ the equivalence gives $N=(N-1)/N$, a well-known exponent.

For the  fundamental solution of the fractional Laplacian diffusion, the decay rate for fixed time is $P_t(x)=O(|x|^{-(N+2s)})$. This gives the formal equivalence between spatial decay rates: $N+2s\sim 2/(1-m)$. It does not agree with the heat equation value in the limit: $m=1$, $s=1$.

One may wonder if we can get complete agreement of exponents and profile functions. This happens in dimension $N=1$ for $s=1/2$ and $m=0$, a very exceptional case. The profile function is the Cauchy distribution
\begin{equation}\label{Cauchy.dis}
P(x)=\frac{1}{\pi}\frac{t}{t^2+|x|^2}\,.
\end{equation}

\subsection{Subcritical, logarithmic and very singular fast diffusion}

The situation becomes much more involved once we cross the value $m_c=(N-2)/N$ for $N\ge 3$. In the range $m_c>m>0$ many of the above  properties do not hold. Thus, there must be a condition on the initial data to guarantee that solutions becomes bounded, hence smooth, for all positive times. If $u_0\in L^P_{loc}(\ren)$ then we need $p>N(1-m)/2$, and the bound is sharp, otherwise regularization need not be true. When $u_0\in L^1(\ren)$, conservation of mass never holds, and in fact such solutions disappear in finite time, a phenomenon called extinction, that is discussed at length in our \cite{VSmoothing}. Obtaining valid versions of the Harnack inequality was challenging in this range \cite{BonVa05, BonVa10}, see also the monograph \cite{DiBGVbk}.

\smallskip

\noindent $\bullet$ {\sc Logarithmic diffusion.}  The limit $m\to 0$ makes sense if we slightly modify the constants in the equation and write it as
\begin{equation}\label{fde.m}
u_t=\nabla\cdot (u^{m-1}\nabla u)=(1/m)\Delta u^m.
\end{equation}
It is proved that the solutions $u_m(x,t)$ with fixed initial data, say bounded, converge as $m\to 0$ to a solution of the logarithmic diffusion equation
\begin{equation}
u_t=\nabla\cdot (\nabla u/ u)=\Delta \log(u)\,,
\end{equation}
famous in 2D as a model for the evolution of the conformal matric by Ricci flow, as proposed by Hamilton \cite{Ham88} in 1988, where $u$ is the conformal factor. A detailed study of the surprising mathematical theory is done in \cite{VSmoothing}, where references are given. The following facts are remarkable:  finite mass solutions are not uniquely determined by the initial data and moreover, they all lose at least $4\pi$ units of mass (which here means surface) per unit time. A very beautiful solution happens when we choose surface loss equal to $8\pi$ and the formula is
\begin{equation}
U(x,t)=\frac{8a(T-t)}{(a^2+|x|^2)^2}, \quad \mbox{with } \ a>0.
\end{equation}
In the geometrical interpretation it describes the shrinking of a perfect 2D ball to a point in time $T>0$. The ball is represented on the plane by stereographical projection.The solution qualifies as another beautiful diffusive pattern, this time it portrays extinction by Ricci flow. We ask the reader to note the difference with the Barenblatt FDE solutions, or with the Cauchy distribution \eqref{Cauchy.dis}.

Note that there is another natural limit as $m\to 0$, namely the equation $u_t=\Delta\,\mbox{\rm sign} (u)$. Though it means no flow for positive data, it has an interesting interpretation in terms of total variation flow for signed data, see \cite{BF2012} in 1D and compare with \cite{RoVa2002}. Total variation flow is an very important subject in itself, related to the $p$-Laplacian, cf. the monograph \cite{ACM2004}.

\noindent $\bullet$  {\sc Super-fast diffusion.} Once we cover $m=0$ the natural question is, can we cover $m<0$. Actually, formula \eqref{fde.m} makes perfect sense and a theory can be developed that extends much of what we have seen in the subcritical case $0<m<m_c$. Some surprises arise in the form of nonexistence for integrable data (which seem in principle the most natural), see \cite{Vnonex92}. The range is called very singular diffusion of super-fast diffusion.

\noindent $\bullet$  {\sc Subcritical asymptotic stabilization.} The absence of Barenblatt solutions makes one wonder what happens for large times in the subcritical, logarithmic and very singular cases. This is a complicated topic, that needs lots of mathematics. We refer to \cite{King93} and \cite{VSmoothing} for the earlier extinction analysis for so-called small solutions, and to \cite{BDGVProc}, \cite{BBDGV2} for stabilization of solution with certain fat tails to the self-similar solutions called pseudo-Barenblatt solutions. The proofs are based on entropy-entropy dissipation methods. References to abundant related work are found.

\subsection{Comments and extensions }
The subject FDE in the lower $m$ ranges is quite rich. Part of the very interesting results concern problems posed in bounded domains, where the discussion is quite different. However, lack of space leads us not to continue the study of the Fast Diffusion range, and we refer to reader to monographs like \cite{DaskKenig, VSmoothing}. But let us just point out that there is no unity in the mathematics of the Fast Diffusion comparable to the Porous Medium range, and a number of critical exponents $m<1$ appear. This is the source of many interesting functional developments and physical phenomena that researchers are still trying to understand.

The contents of the two last sections on PME and FDE can be translated to a large extent to the study of the evolution $p$-Laplacian equation (PLE), though some remarkable differences exist. We refer to the  book \cite{VSmoothing} for an account of our ideas. There is even a transformation that maps all radial solutions of the PME to the corresponding class of the PLE, see \cite{ISV08}. Of course, $p=m+1$ by dimensional considerations, but the transformation changes also the space dimension. If $N>2$, then the corresponding PLE dimension is $N'=(N-2)(m+1)/2m$.

Some  studies deal with the Doubly Nonlinear Equation $u_t=\nabla\cdot(|D(u^m)|^{p-2} D(u^m))$.
See \cite{StanVaDNL} for a recent work.


\section{\textsf{Nonlinear fractional diffusion. Potential model}}
\label{sec.pmfp}

The combination of fractional diffusion and porous medium nonlinearities gives rise to interesting mathematical models that have been studied in the last decade both because of a number of scientific applications and for their mathematical properties. Two main models will be discussed below; a mechanical model has been developed in collaboration with Luis Caffarelli in Texas, and can be called {\sl porous medium flow with fractional potential pressure} (or more generally, with nonlocal pressure); it has surprising properties. The other one has been developed later but it has better analytical properties. For convenience we will call them here {\bf PMFP} and {\bf FPME}.\footnote{In \cite{Vaz14} they were called Type I and Type II in reverse order.} We will also examine models that interpolate between both.

\subsection{Porous medium diffusion with nonlocal pressure}

We devote this section to introduce model {\bf PMFP}.  It arises from the consideration of a continuum, say, a fluid, represented by a  {\sl density} distribution $u(x,t)\ge 0$ that
evolves with time following a {\sl velocity field } $\bf v(x,t)$,
according to the continuity equation
$$
u_t+\nabla\cdot(u\, {\bf v})=0.
$$
We assume next that $\bf v$ derives from a potential, ${\bf
v}=-\nabla p$, as happens in fluids in porous media
according to Darcy's law, and in that case $p$ is the {pressure}. But
potential velocity fields are found in many other instances, like
Hele-Shaw cells, and other recent examples.

 We still need a closure relation to relate  $p$ to $u$. In the
case of gases in porous media, as modeled by Leibenzon and Muskat,
the closure relation takes the form of a state law { $p=f(u)$,} where
$f$ is a nondecreasing scalar function, which is linear when the
flow is isothermal, and a power of $u$ if it is adiabatic. The PME follows.
The linear relationship happens also in the simplified description of water infiltration in an almost horizontal soil layer according to Boussinesq's modelling. In that case we get the standard porous medium equation, $u_t=c\Delta (u^2)$. See Section \ref{sec.pme} on the PME or \cite{Vazpme} for these and many other applications.

The diffusion model with nonlocal effects proposed in 2007 with Luis Caffarelli  uses
the first steps of the derivation of the PME, but it differs by using a closure relation of the form \  $p={\cal K}(u),$  where $\cal K$ is a linear integral operator, which we assume in practice to be the inverse of a fractional Laplacian. Hence, $p$ es related to $u$ through a fractional potential operator,
${\cal K}=(-\Delta)^{-s}$, $0<s<1,$ with kernel
$$
k(x,y)= c|x-y|^{-(n-2s)}
$$
(i.e., a Riesz operator). We have $(-\Delta)^s p = u$.
This introduces long-distance effects in the model through the pressure, and we end up with a nonlocal model,  given by the system
\begin{equation}\label{eq.pmp}
\fbox{\large  $ u_t=\nabla\cdot(u\,\nabla p), \quad p={\cal K}(u)$ }
\end{equation}
where $u$ is a function of the variables $(x,t)$ to be thought of
as a density or concentration, and therefore nonnegative, while
$p$ is the nonlocal pressure, which is related to $u$ via a linear operator
$\cal K$. We can write:  $u_t=\nabla\cdot(u\,\nabla (-\Delta)^{-s}u)$. A technical observation: there are problems in defining  $(-\Delta)^{-s}u$ in 1D since the kernel may be too singular, but then $\nabla (-\Delta)^{-s}u$ is always well defined, which is enough to perform the calculations that will be commented upon below.

The problem is posed for $x\in \ren$, $n\ge 1$, and $t>0$, and we
give initial conditions
\begin{equation}\label{eq.ic}
u(x,0)=u_0(x), \quad x\in \ren,
\end{equation}
where $u_0$ is a nonnegative, bounded and integrable function in $\ren.$

\subsection{Applied motivation and variants} $\bullet$  {\sc Particle systems with long range
interactions}. Equations of the more general form
$$
u_t=\nabla \cdot (\sigma(u)\nabla {\cal L}u)
 $$
have appeared in a number of applications to the macroscopic evolution of particle systems. Thus, Giacomin and Lebowitz \cite{GL97}, 1997, consider  a lattice gas with general short-range interactions and a Kac potential, and passing to the limit, the macroscopic density profile
$\rho(r,t)$ satisfies the equation
\begin{equation}
\qquad \frac{\partial \rho}{\partial t}=\nabla   \cdot\left[\sigma_s(\rho)\nabla
\frac{\delta F(\rho)}{\delta \rho }\right]
\end{equation}
where $\sigma_s(\rho)$ may be degenerate. See also \cite{GLM2000}.

$\bullet$ {\sc Modeling dislocation dynamics as a continuum.} Following old modeling by A. K. Head in \cite{Head}, Biler-Karch-Monneau \cite{BKM10} considered the one-dimensional case of  model \eqref{eq.pmp}. By integration in $x$ they introduced viscosity solutions \`a la Crandall-Evans-Lions. They prove that uniqueness holds, which is very satisfying property. But the corresponding mathematical model in several space dimensions looks quite different from \eqref{eq.pmp}.

\smallskip

$\bullet$  {\sc Hydrodynamic limit for $s=1$.} This is a very interesting limit case. Putting $s=1$ makes us lose the parabolic character of the flow that becomes hyperbolic. In 1D the situation is rather trivial since when we put $p=(-\Delta)^{-1}u$ we get $p_{xx}=-u_x$, and then
 $$
 u_t  = (u\, p_x)_x= u_x\,p_x- u^2
 $$
Moreover, if $v=-p_x=\int u\,dx$, we have
$$
v_t=up_x+ c(t)=-v_x v+c(t),
$$
For $c=0$ this is the  Burgers equation $v_t+vv_x=0$, which generates shocks
in finite time but only if we allow for $u$ to have two signs.

In several dimensions the issue becomes much more interesting because it does not reduce to a simple Burgers equation. We have
\begin{equation}
 u_t  = \nabla\cdot(u\,\nabla p)=\nabla u\cdot\nabla p- u^2, \qquad p=(-\Delta)^{-1}u\,,
\end{equation}
A  very close version to this model has appeared in superconductivity (the Chapman-Rubinstein-Schatzman-E model) see  \cite{LinZh00, W-E}, and Ambrosio-Serfaty \cite{AmSr08}. In that application $ u$ describes the vortex density. Gradient flow structure for this example is established in  \cite{AMS}.

\smallskip

$\bullet$ {\sc The PME limit.} If we take $s=0$, ${\cal K}=$ the
identity operator, we get the standard porous medium equation,
whose behaviour is well-known. Therefore we can see the PMFP equation as
a nonlinear interpolation between the PME and the hydrodynamic limit, $s=1$.

\smallskip

$\bullet$  More generally, it could be assumed that ${\cal K}$ is an operator of
integral type defined by convolution on all of $\RR^n$, with the
assumptions that is positive and symmetric. The fact the ${\cal K}$ is a
homogeneous operator of degree $2s$, $0<s<1$, will be important in
the proofs.  An interesting variant would be the Bessel kernel \
 ${\cal K}=(-\Delta+cI)^{-s}$. We are not exploring such extensions.

\subsection{Mathematical results}

Early results on the PMFP have been reported in Proceedings from the Abel Symposium \cite{Vaz12Abel}, and then in \cite{Vaz14}, so we will concentrate on general facts and only develop in more detail some of the new material. For applications of nonlinear nonlocal diffusion see also \cite{Caff12Abelsurv}.

$\bullet$ In paper \cite{CV1} Luis Caffarelli and the author established the  existence of weak energy solutions, the basic properties of the solutions, like  conservation of mass
\begin{equation}
\frac{d}{dt}\int u(x,t)\,dx=0\,,
\end{equation}
the two energy estimates
\begin{equation}
\frac{d}{dt}\int u(x,t)\log u(x,t)\,dx=-\int |\nabla Hu|^2\,dx\,,
\end{equation}
where $ H=(-\Delta)^{-s/2} $, and
\begin{equation}
\frac{d}{dt}\int |Hu(x,t)|^2\,dx=-2\int u|\nabla Ku|^2\,dx, \quad K=(-\Delta)^{-s}.
\end{equation}
A number of usual properties in diffusive processes do hold here like conservation of positivity, as well as $L^p$ decay. But we also found lack of a general comparison principle, a major difficulty in developing the theory (such a drawback will not be shared by the second model, FPME). And we could not prove uniqueness for general solutions in several space dimensions.

A main goal in the study of this model was to determine whether or not the {\sl property of finite propagation } holds. The answer turned out to be yes. This is not clear in principle due to the competition between the slow propagation of the PME part with the infinite propagation of the fractional operator (amounting to long distance effects). The lack of plain comparison made the proof difficult, and the difficulty was surmounted by a novel use of the methods of viscosity solutions. Summing up, the degenerate character of the PME wins. On the contrary, infinite propagation was later proved to be true for FPME.

$\bullet$ In a second contribution \cite{CV2} we explored the long-time behaviour in two steps. We first established the  existence of self-similar profiles, so-called Fractional Barenblatt solutions
$$
U(x,t)=t^{-\alpha}F(x\,t^{-\beta}), \qquad \beta=\frac{1}{N+2-2s}, \ \alpha=N\beta,
$$
The profile $F$ is compactly supported, a clue to the finite propagation property, and  is the solution of a certain fractional obstacle problem.  A different proof in dimension 1 follows from paper \cite{BKM10}. The authors of \cite{BIK11} found the self-similar Barenblatt profiles in all dimensions with explicit formulas:
$F(x)=(A-B\,|x|^2)_+^{1-s}$.

Then we introduced the renormalized Fokker-Planck equation and used a suitable entropy functional and  proved stabilization of general solutions to the  previous profiles that we called fractional Barenblatt profiles. All this is carefully explained in  \cite{Vaz12Abel}.

$\bullet$  The next issue in the programme was the regularity of the solutions. It was studied in a paper with L. Caffarelli and F. Soria, \cite{CSV}. Proving boundedness for solutions with integrable data in $L^p$, $1\le p\le \infty$ was  an important step in the this theory. We can dispense with the extension method for fractional Laplacians by using energy estimates based on the properties of the quadratic and bilinear forms associated to the fractional operator, and then the iteration technique.

 {\bf Theorem}  {\sl Let $u$ be a weak solution the initial-value problem for the  PMFP  with data $u_0\in L^1(\mathbb{R}^n)\cap L^\infty(\mathbb{R}^n)$,  as constructed before. Then, there exists a positive constant $C$ such that for every $t>0$
\begin{equation}\label{form:L-inf}
\sup_{x\in\mathbb{R}^n}|u(x,t)|\le C\,t^{-\alpha }\|u_0\|_{L^1(\mathbb{R}^N)}^{\gamma}
\end{equation}
with $\alpha=N/(N+2-2s)$, $\gamma=(2-2s)/(N+2-2s)$. The constant $C$ depends only on $N$ and $s$.}

The major step is then proving  $C^\alpha$ regularity.  The proof uses the DeGiorgi method with careful truncations together with very sophisticated energy methods that have to overcome the difficulties of both nonlinearity and nonlocality. A number of ideas come from  Caffarelli-Vasseur \cite{CaffVass} and \cite{CafChanVas} with difficult modifications due to the degenerate nonlinearity. The theory can be extended to data $u_0\in L^1(\ren)$, $u_0\ge 0$, giving global existence of bounded weak solutions.

\subsection{Energy, bilinear forms and fractional Sobolev spaces}
The previous results are obtained in the framework of weak energy solutions:  The basis of the boundedness analysis is a property that goes beyond the definition of weak solution. The general energy property is as follows: for any real smooth function $F$ and such that $f=F'$ is bounded and nonnegative, we have   for every $0\le t_1\le t_2\le T$,
\begin{equation*}\label{entro}
\begin{array}{ll}
\int F(u(t_2))\, dx -\int F(u(t_1))\, dx & = - \int_{t_1}^{t_2}\int \nabla [f(u)] u \nabla p\, dx\, dt= \\
& -\int_{t_1}^{t_2}\int \nabla h(u) \nabla (-\Delta)^{-s} u\, dx\, dt\,,
\end{array}
\end{equation*}
where  $h$ is a function  satisfying  $h'(u)= u\,f'(u)$. We can write the last integral in terms of a  bilinear form
\begin{equation*}
\int \nabla h(u) \nabla (-\Delta)^{-s} u\, dx = \mathcal{B}_s (h(u), u)
\end{equation*}
  This bilinear form $\mathcal{B}_s$ is defined  as
\begin{equation*}
\begin{array}{l}
\mathcal{B}_s(v,w)= C \displaystyle \iint  \nabla v(x)\frac{1}{|x-y|^{N-2s}} \nabla w(y)\, dx \, dy = \\
\displaystyle \iint {\cal N}_{-s} (x,y) \nabla v(x) \nabla w(y) \, dx\, dy
\end{array}
\end{equation*}
where ${\cal N}_{-s}(x,y)= C\,|x-y|^{-(N-2s)}$ is the kernel of operator $(-\Delta)^{-s}$.
 After some integrations by parts we also have
\begin{equation}\label{BB}
\mathcal{B}_s (v,w)= C_{n,1-s} \iint(v(x)-v(y)) \frac{1}{|x-y|^{n+ 2(1-s)}} (w(x)-w(y))\, dx\, dy\end{equation}
since $-\Delta {\cal N}_{-s}={\cal N}_{1-s}$.  It is well known  that $\mathcal{B}_s (u,u)$ is an equivalent norm for the fractional Sobolev space $W^{1-s,2}(\ren)$. This is the way the fractional Sobolev spaces appear, as dissipated energies that will guarantee compactness in the arguments, see \cite{CSV}. Fractional Sobolev spaces with a view to their use in PDEs are discussed in \cite{DNPV11}.

\subsection{More recent work }

$\bullet$ The particular value $s=1/2$ of the fractional exponent turned out to be extremely delicate for the regularity analysis of \cite{CSV} and needed a further article with new geometrical ideas, \cite{CV2}. Briefly stated, there are some terms in the energy estimates that come from the tails of the solutions and cannot be suitably controlled in an iterative way. An iterated geometrical coordinate transformation  allows to eliminate them at the cost of a {\sl distorted geometry}.

 \smallskip

 $\bullet$  The Hydrodynamic Limit $s\to 1$ was studied by  S. Serfaty and the author in  \cite{SerVaz}.
 We pass to the limit and construct a theory of existence, uniqueness and estimates for the hydrodynamic limit problem. It is interesting to note that the asymptotic attractor is a selfsimilar vortex of the form
 $$
 U(x,t)= t^{-1} F(x/t^{1/N}),
 $$
 and $F$ is the characteristic function of a ball. Therefore, even continuity is lost in the regularity of the solutions. This is not a contradiction since the limit equation is no longer parabolic. Our work is related to work on aggregation models by Bertozzi et al.  \cite{BLL}.

 \smallskip

$\bullet$   We posed the question of possible rates in the asymptotic convergence to  selfsimilar solutions of Barenblatt type of papers \cite{CV2, BKM10}. This question was partially solved in a paper  \cite{CH2013} with  J. A. Carrillo, Y. Huang and M. C. Santos, where we showed exponential convergence towards stationary states for the Porous Medium Equation with Fractional Pressure in 1D. The many-dimensional case seems to be a difficult open problem,  it is tied to some functional inequalities that are not known. Our analytical approach does not seem to apply either.


\subsection{Additional work, open problems}

\quad $\bullet$ The questions of uniqueness and comparison are solved in dimension $N=1$ thanks to the trick of integration in space used by Biler, Karch, and Monneau \cite{BKM10}. New tools are needed to make progress in several dimensions.

Recent uniqueness results are due to Zhou, Xiao, and Chen, \cite{ZXC}. They obtain local in time strong solutions in Besov spaces. Thus, for initial data in $B^\alpha_{1,\infty}$ if $1/2\le s<1$ and $\alpha>N+1$ and $N\ge 2$. Therefore, Besov regularity implies uniqueness for small times.

$\bullet$  The fractional Burgers connection was explored in \cite{CasCor08} for $N=1$, $s=1/2$ where $\partial_x(-\Delta)^{-1/2}=-H$, the Hilbert transform.

$\bullet$ The study of the free  boundary is in progress, but regularity  is still open for small $s>0$.

$\bullet$   The gradient flow structure of the {\bf PMFM} flow in Wasserstein metrics has been recently established  by S. Lisini, E. Mainini and A. Segatti in \cite{LMS16}. For the general approach see the monograph \cite{AGS}. Previous work in 1D was due to by J. A. Carrillo et al.

$\bullet$ The problem in a bounded domain with Dirichlet or Neumann data has not been studied, to our knowledge.

$\bullet$ Good numerical algorithms and studies are needed.

\subsection{Elliptic  nonlinear nonlocal models.} The interest in using { fractional Laplacians} in
modeling diffusive processes has a wide literature,
especially when one wants to model long-range diffusive interactions,
and this interest has been activated by the recent progress in the
mathematical theory, in the form of a large number works on elliptic equations,
mainly of the linear or semilinear type, as well as free boundary problems, like obstacle problems. There are so many works on the subject that we cannot refer them here. Let us mention the survey paper \cite{Caff12Abelsurv} by L. Caffarelli, that contains  a discussion of the properties of solutions to several non-linear elliptic equations involving diffusive processes of non-local nature, including reference to drifts and game theory.

\section{\textsf{The FPME model and the mixed models}}
\label{sec.mod2}

Another natural model for the combination of fractional diffusion and porous medium nonlinearities is the equation that we will call fractional porous medium equation:
 $\partial_t u+(-\Delta)^s (u^m)=0$. In order to be mathematically precise we write the equation as
\begin{equation}
\fbox{ $ u_t+(-\Delta)^{s}(|u|^{m-1}u)=0 \ $}
\end{equation}
with $0<m<\infty$ and  $0<s<1$. We will take initial data in $u_0\in L^1(\ren)$ unless mention to the contrary. Normally, $u_0, u\ge 0$. We will refer to this model as  {\bf FPME} for easy reference in this paper.
Mathematically, it looks a more direct generalization of the linear fractional heat equation than the potential model {\bf PMFP} studied in the previous section.

This model represents another type of nonlinear interpolation with parameter $s\in (0,1)$, this time between the PME \ $u_t -\Delta(|u|^{m-1}u)=0$ for $s=1$ and the plain absorption ODE \  $u_t+|u|^{m-1}u=0$ for $s=0$.

We have written a detailed description of this model in the survey paper \cite{Vaz14}, where we give references to the physical motivations, among them \cite{AC09, JKOlla, Jara1, Jara2}, the literature, and the mathematical developments until 2013 approximately. See also Appendix B of \cite{BonVa10}. Therefore,  we will mention the main items of the research, the references and general ideas, and then proceed to give notice of recent work, that covers different directions.

\subsection{Mathematical theory of the FPME}

A complete analysis of the Cauchy problem posed for $x\in\RR^N$, $t>0$, with initial data in $L^1(\ren)$ was done in two very complete papers coauthored with A. de Pablo, F. Quir\'{o}s,  and A. Rodr\'{\i}guez: \cite{pqrv1} in (2011) and  \cite{pqrv2} (2012). Using the Caffarelli-Silvestre extension method and the B\'enilan-Brezis-Crandall functional semigroup approach, a weak energy solution is constructed, and $u\in C([0,\infty):L^1(\ren))$. Moreover, the set of solutions forms a semigroup of ordered contractions in $L^1(\ren)$. This is the first instance of a `better behaviour' than model {\bf PMFP}.

\smallskip

$\bullet$  The second big difference is that {\sl Nonnegative solutions have infinite speed of propagation for all $m$ and $s$, so that there is no nonnegative solution with compact support
(we mean, in the space variable).} Actually, a  very important property of Model PMFP with Caffarelli is that solutions with compactly supported initial data  do have the compact support property (i.e., they stay compactly supported for all times).

Even is propagation is always infinite in this model, we still use the name `fractional FD' for the range $m<1$ because of the many analogies with the usual FDE.

\smallskip

$\bullet$ On the other hand, some  properties are similar in both models:  Conservation of mass holds for all $m\ge 1$, and even for some $m<1$ close to 1; the $L^1-L^\infty$ smoothing effect works; and the $C^\alpha$ regularity holds also (unless $m$ is  near 0 and solutions are not bounded). Comparison of the models {\bf PMFP} and {\bf FPME}  is quite interesting and has been pursued at all levels.

\smallskip

$\bullet$  The question of existence of classical solutions and higher regularity for the FPME and the more general model
$$
\partial_t u +(-\Delta)^s\Phi(u)=0
$$
(where $\Phi$ is a monotone real function with $\Phi'\ge 0$) has been studied in two papers with the same authors (A.deP., F.Q., A.R., J.L.V.). The first paper, \cite{pqrv3}, treats the model case $\Phi(u)=\log(1+u)$, which is interesting as a case of log-diffusion. The second treats general nonlinearities $\Phi$ and proves higher regularity for nonnegative solutions of this fractional porous medium equation, \cite{vpqr}. This is a very delicate result. There is an extension of this result to prove $C^\infty$ regularity to solutions in bounded domains by  Bonforte, Figalli and Ros-Ot\'on, \cite{BFR2017}.

\smallskip

$\bullet$ Our paper \cite{VazBar2012} deals with the construction of  what we call  the fractional Barenblatt solution of the FPME, which has the also self-similar form:
\begin{equation}
 U(x,t)=t^{-\alpha}F(xt^{-\beta})
\end{equation}
The construction works for $m>m_c=(N-2s)/N$, a range that is optimal that reminds us of the Fast Diffusion Equation, Section \ref{subsec.fde}.  The difficulty is to find $F$ as the solution of an elliptic nonlinear equation of fractional type. Such profile is not explicit as in the {\bf PMFP} model (it is only for some very special exponents \cite{Huang}). In any case, $F$ has behaviour like a power tail
 $$
 F(r)\sim r^{-(N+2s)}\,.
 $$
This is important for the applications and it the same as the one predicted by Blumental for the linear fractional kernel. This asymptotic spatial behaviour holds for all $m\ge 1$, and even for some $m<1$,  but  not for some  fast diffusion exponents $m_c<m<1$ (see what happens then in  \cite{VazBar2012}).

 Asymptotic behaviour as $t\to\infty$ follows, and this Barenblatt pattern is proved to be an attractor, as we were expecting from what has been seen along this whole text. The result holds for $m>m_c$. Open problem: Rates of  convergence have not been found, and this is an interesting open problem.

 Extinction in  finite time is proved for exponents $0<m<m_c$. The corresponding stabilization process must be studied.

\smallskip

$\bullet$  Another direction concerns  regularity at the local or global level. In collaboration with M. Bonforte we have obtained a priori upper and lower estimates of intrinsic, local type for this problem posed in $\ren$,  \cite{BV2012}. Quantitative positivity and Harnack Inequalities follow.  Against some prejudice due to the nonlocal character of the diffusion, we are able to obtain them here for fractional PME/FDE using a technique of weighted integrals  to control the tails of the integrals in a uniform way. The novelty are the weighted functional inequalities. This also leads to existence of solutions in weighted $L^1$-space
for the fast diffusion version FPME, a restriction that does not appear in the standard FDE.

More recent, very interesting work on bounded domains is reported in Section \ref{sec.bdd.domain}.

\smallskip

$\bullet$  Symmetrization (Schwarz and Steiner types). This is a project  with  B. Volzone \cite{VazVol, VazVol2}. Applying usual symmetrization techniques is not easy and we have found a number of open problems. It turns out that Steiner symetrization works and it does much better for fractional FDE than for the fractional PME range. This work was followed by recent collaboration with Y. Sire and B. Volzone to apply the techniques to the fractional Faber-Krahn inequality, \cite{SVV}.

\smallskip

$\bullet$ We have also investigated very degenerate nonlinearities, like the {\sl Mesa Problem}. This is the limit of FPME with $m\to\infty$.  We have studied this limit  in  \cite{Vaz15}, and the limit flow characterized by the solution of a fractional obstacle problem, that is related to the obstacle problem for the PMFP that was described in \cite{CV2}.

\smallskip

$\bullet$  Numerics for the nonlinear nonlocal diffusion models is being done by a number of authors at this moment by :  Nochetto et al. \cite{NOS13, NOS16}, Teso \cite{Teso14}.

\smallskip

$\bullet$  Fast diffusion and extinction. Very singular fast diffusion.
Paper with Bonforte and Segatti  \cite{BSgV2016}, on non-existence due to instantaneous extinction, which is the common rule in very singular fractional fast diffusion as shown for standard diffusion in \cite{Vnonex92}. Paper \cite{Vaz16jee} shows the existence of maximal solutions for some very singular nonlinear fractional diffusion equations in 1D in some borderline cases, this is an exception.

\smallskip

$\bullet$ We have looked at the phenomenon of KPP propagation in linear and nonlinear fractional diffusion with the particular reaction proposed by Kolmogorov-Piskunov-Petrovskii and Fisher (1938). In the case of standard linear diffusion travelling waves appear and serve as asymptotic attractors. Cabr\'e and Roquejoffre \cite{CabreRoquejoffre1, CabreRoquejoffreArxiv} studied the diffusion equation with linear fractional diffusion and  KPP reaction and showed that there is no traveling wave propagation, and in fact the level sets move out at an exponential rate for large times. The results are extended to nonlinear fractional diffusion of the FPME type for all values of the exponents in work with Stan, \cite{StanVazquezKPP}.

\subsection{Mixed Models}\label{sde.mixed}

$\bullet$ The potential model {\bf PMFP} given in \eqref{eq.pmp} is generalized into {\bf PMFP'}
\begin{equation}\label{m1ext}
u_t=\nabla \cdot (u^{m-1}\nabla (-\Delta)^{-s} u)
\end{equation}
 with $m>1$. This is an extension that accepts a general exponent $m$, so that the comparison of both models may take place on more equal terms.

 The most interesting question seems to be deciding if there is  finite and infinite propagation for {\bf PMFP'}. Recent works with D. Stan and F. del Teso \cite{StTsVz.CRAS} and  \cite{StTsVz16} show that finite propagation is true for $m\ge 2$ and propagation is infinite is $m<2$. This is quite different from the standard porous medium case $s=0$, where $m=1$ is the dividing value of the exponent as regards propagation. The problem with existence is delicate for large $m$ and is treated in  a further paper \cite{StTsVz17}.

An interesting and unexpected aspect of the theory is the existence of  a  transformation
that maps self-similar solutions of the FPME with $m\ge 1$ into solutions of the same type for
model {\bf PMFP'} with exponent $1<m<2$. This applies in particular to the Barenblatt solutions constructed in \cite{VazBar2012}. The transformation is established in \cite {StTsVz15} and is quite useful in showing that {\bf PMFP'} has infinite propagation in that range of parameters.

\medskip

 $\bullet$  Work by  Biler-Imbert-Karch \cite{BIK2015} deals with the variant
\begin{equation}
u_t=\nabla \cdot (u\nabla (-\Delta)^{-s} u^{m-1})\,.
\end{equation}
They construct a family of nonnegative explicit compactly supported self-similar solutions which are a generalization of the well-known Barenblatt profiles for the classical porous-medium equation. They also establish the existence of sign-changing weak solutions to the initial-value problem, which satisfy sharp hypercontractivity $L^1$-$L^p$ estimates.

\medskip

 $\bullet$  Reference \cite{StTsVz15} also treats  on the double exponent model
\begin{equation}\label{doble}
\partial_t u +\nabla(u^{m-1}\nabla (-\Delta)^{-s}u^{n-1})=0
\end{equation}
that generalizes all the previous models. The paper discusses self-similar transformations and finite propagation. The transformation of self-similar solutions indicates that finite propagation holds for $m\ge 2$, while $n>1$ does not count. The graphic of parameters in \cite{StTsVz15} gives a very clear scheme of these transformations.

\subsection{Some related directions }

There is work on equations with other nonlocal linear operators, and also on
equations with lower order terms, leading to reaction-diffusion and blow-up.
Nonlinear diffusion and convection is treated in \cite{CiJak11, ACJ12}.
The chemotaxis systems have been studied with nonlocal and/or nonlinear diffusion, like \cite{Esc06, BilerWu, BourCal2010, LiRoZh2010}. We also have geometrical flows, like the fractional Yamabe problem (to be mentioned below). And there are a number of other options.

\section{\textsf{Operators and Equations in Bounded Domains}}
\label{sec.bdd.domain}

 We have presented  different definitions of the fractional Laplacian operator acting in $\ren$ in Section \ref{sec.fraclap}, and we have mentioned that all these versions are equivalent. However, when we want to pose a similar operator in a bounded domain $\Omega\subset \ren$ we have to re-examine the issue, and several non-equivalent options appear. This enlarges the theory of evolution equations of fractional type on bounded domains, and the recent literature has taken it into account. Actually, there is much recent progress in this topic and the next second subsection will describe our recent contributions. A large class of related nonlocal diffusive operators can be considered in the same framework.

 \subsection{The linear operators} There are a number of definitions that have been suggested for the fractional Laplacian operator (FLO) acting on a bounded domain $\Omega$. The ones we consider here are naturally motivated, and they give rise to different  operators. We will mention three basic options, two of them are mostly used.

 \medskip

\noindent {\bf  \large  The  Restricted Fractional Laplacian operator (RFL)}.
It is the simplest option. It acts on functions $g(x)$ defined in $\Omega$ and extended by zero to the complement, and then the whole  hypersingular integral of the Euclidean case is used. Therefore, it is just the fractional Laplacian in the whole space ``restricted'' to functions that are zero outside $\Omega$.
\begin{equation}
(-\Delta_{|\Omega})^{s}  g(x)= c_{N,s}\mbox{
P.V.}\int_{\mathbb{R}^N} \frac{g(x)-g(z)}{|x-z|^{N+2s}}\,\rd z\,,\qquad\mbox{with }\mbox{supp}(g)\subset \overline{\Omega}\,.
\end{equation}
Here, $s\in(0,1)$ and $c_{N,s}>0$ is a normalization constant. It is shown that, thus defined,  $(-\Delta_{|\Omega})^s$ is a self-adjoint operator on $L^2(\Omega)$ with a discrete spectrum, with eigenvalues
$$
0<\overline{\lambda}_1\le\overline{\lambda}_2\le \ldots\le\overline{\lambda}_j\le\overline{\lambda}_{j+1}\le \ldots,
$$
satisfying $\overline{\lambda}_j\asymp j^{2s/N},$ for $j\gg 1$. The corresponding eigenfunctions $\overline{\phi}_j$ are only H\"older continuous up to the boundary, namely $\overline{\phi}_j\in C^s(\overline{\Omega})$\,,  \cite{RosSer}.

An important issue is the way in which the additional conditions (formerly boundary conditions) are implemented for the RFL. It usually takes the form of exterior conditions:
\begin{equation}
u(t,x)=0\,,\; \qquad \mbox{in }(0,\infty)\times \big(\RR^N\setminus\Omega\big)\,.
\end{equation}

The behavior of the Green function $G$ plays an important role in the corresponding PDE theory. It satisfies a strong behaviour condition,  that we call {(K4) condition}:
\begin{equation}\tag{\rm K4}
G(x,y)\asymp \frac{1}{|x-y|^{N-2s}}
\left(\frac{\p(x)}{|x-y|^\gamma}\wedge 1\right)
\left(\frac{\p(y)}{|x-y|^\gamma}\wedge 1\right)\,,
\end{equation}
where $\delta(x)$ is the distance from $x\in\Omega$ to the boundary. The exponent $\gamma$ will play a role in the results derived from the kernel. In the RFL we have $\gamma=s$.

\noindent\textbf{References.}  There is an extensive literature on the RFL operator and the corresponding $\alpha$-stable process in the probability literature. The interested reader is referred \cite{BSV2015}where we have commented on relevant works in that direction. 
 \normalcolor

\medskip

\noindent {\bf \large  The  Spectral Fractional Laplacian operator (SFL).}
It is defined by the two equivalent expressions
\begin{equation}\label{sLapl.Omega.Spectral}
\displaystyle(-\Delta_{\Omega})^{s}
g(x)=\sum_{j=1}^{\infty}\lambda_j^s\, \hat{g}_j\, \phi_j(x)
=\frac1{\Gamma(-s)}\int_0^\infty
\left(e^{t\Delta_{\Omega}}g(x)-g(x)\right)\frac{dt}{t^{1+s}}\,,
\end{equation}
where $\Delta_{\Omega}$ is the classical Dirichlet Laplacian on the domain $\Omega$,
and $\hat{g}_j$ are the Fourier coefficients of $f$
\[
\hat{g}_j=\int_\Omega g(x)\phi_j(x)\dx\,,\qquad\mbox{with}\qquad \|\phi_j\|_{L^2(\Omega)}=1\,.
\]
In this case the eigenfunctions $\phi_j$ are  the same as in the Dirichlet Laplacian, smooth as the boundary of $\Omega$ allows. Namely, when $\partial\Omega$ is $C^k$, then  $\phi_j\in C^{\infty}(\Omega)\cap C^k(\overline{\Omega})$ for all $k\in \mathbb{N}$\,.
The eigenvalues are powers $\lambda_j^s$ of the standard eigenvalues $0<\lambda_1\le\lambda_2\le \ldots\le\lambda_j\le\lambda_{j+1}\le \ldots$ and $\lambda_j\asymp j^{2/N}$.
It is proved that the eigenvalues of the RFL are smaller than the ones of SFL: $\overline{\lambda}_j\le
\lambda_j^s $ for all $j\ge 1$, \cite{ChSong2005}.

Lateral boundary conditions for the SFL are different from previous case. They can be read from the   boundary conditions of the Dirichlet Laplacian by the semigroup formula. They are often defined by means of the equivalent formulation that uses the Caffarelli-Silvestre extension defined in a cylinder adapted to the bounded domain, as done in \cite{BCPS, Cabre-Tan, StingaTorrea2010}. If $U$ is the extended function, then we impose $U=0$ on the lateral boundary $x\in \partial \Omega$, $y>0$.

 The Green function of the SFL satisfies the strong assumption (K4), this time with exponent $\gamma=1$.

\medskip

\noindent {\bf  Remarks.} Both  SFL and  RFL admit another possible definition using the so-called Caffarelli-Silvestre extension. They are the two best known options for a FLO. The difference between RFL and SFL seems to have been well-known to probabilists, it was discussed later in PDEs, see  Servadei-Valdinoci \cite{SVal4}, Bonforte and the author  \cite{BVdir2013}, and   Musina-Nazarov \cite{MuNa14}. In this last work the denomination {\sl Navier fractional Laplacian} is used.  The debate about the proper names to be used is not settled.

\medskip

\noindent {\bf \large The Censored Fractional Laplacians (CFL)}.
 This is another option appearing in the probabilistic literature, it has been introduced in 2003 by Bogdan, Burdzy and Chen, \cite{BBChen2003}. The definition is
\begin{equation}
\mathcal{L} g(x)=\mathrm{P.V.}\int_{\Omega}\left(g(x)-g(y)\right)
\frac{a(x,y)}{|x-y|^{N+2s}}\dy\,,\qquad\mbox{with }\frac{1}{2}<s<1\,,
\end{equation}
 where $a(x,y)$ is a measurable, symmetric function bounded between two positive constants, satisfying some further assumptions; for instance $a\in C^{1}(\overline{\Omega}\times\overline{\Omega})$. In the simplest case we put  $a(x,y)=$ constant. On the other hand, \cite{BBChen2003} point out that in the excluded range  $s\in(0,1/2]$ the censored 2s−stable process is conservative and will
never approach the boundary. The CFL is also called {\sl regional fractional Laplacian}.


  \medskip

 The Green function $G(x,y)$ satisfies condition (K4) with $\gamma=s-\frac{1}{2}$, as proven by Chen, Kim and Song \cite{CKS2010b}. See also \cite{GMa}.

\medskip

\noindent {\bf Note.}  We have presented 3 models of Dirichlet fractional Laplacian.
The estimates (K4) show that they are of course \textbf{not equivalent}.
Our work described in the next subsection applies to those operators and a number of other variants, that are listed in \cite{BV.Bdd2016} and \cite{BFV2016}. For instance, sums of operators of the above types and powers of said operators are included.


\subsection{Nonlocal diffusion of porous medium type on bounded domains}
\label{ssec.ndbdd}

We report here on very recent work done in  collaboration with
 M. Bonforte, and also Y. Sire and  A. Figalli, on nonlinear evolution equations of porous medium type posed in bounded domains and involving fractional Laplacians and other nonlocal operators. The papers are  \cite{BSV2015}, \cite{BVdir2013}, \cite{BV.Bdd2016},  \cite{BFR2017}, and  \cite{BFV2016}.

We develop a new programme for  nonlocal porous medium equations on bounded domains aiming at establishing existence, uniqueness, positivity, a priori bounds, regularity, and asymptotic behaviour for a large class of equations of that type in a unified way. We include the set of suitable versions of FLO in a bounded domain. The main equation is written in abstract form as
\begin{equation}
\partial_t u+{\mathcal L} \Phi(u)=0\,,
\end{equation}
where $\Phi$ a continuous and nondecreasing real function, most often  a power function.

\smallskip

\noindent $\bullet$ A problem to be settled first is the suitable concept of solution. We use the ``dual'' formulation of the problem and the concept of {\sl weak dual solution}, introduced in \cite{BVdir2013}, Definition 3.4, which extends the concept of weighted very weak solution used before. In brief, we use the linearity of the operator ${\mathcal L}$ to lift the problem to a problem for the potential function
$$
U(x,t)=\int_\Omega u(y,t)G(x,y)dy
$$
where $G$ is the elliptic Green function for ${\mathcal L}$. Then $\partial_t U=-\Phi(u).$

\smallskip

\noindent $\bullet$ {\sc Class of operators.} In our recent work we have extended the evolution theory to cover a wide class of linear operators $\mathcal{L}$ that satisfy the following conditions. $\mathcal{L}: {\rm dom}(\mathcal{L})\subseteq{\rm L}^1(\Omega)\to{\rm L}^1(\Omega)$ is assumed to be densely defined  and sub-Markovian, more precisely, it satisfies (A1) and (A2):
\begin{enumerate}
\item[(A1)] $\mathcal{L}$ is $m$-accretive on $\rm L^1(\Omega)$;
\item[(A2)] If $0\le f\le 1$ then $0\le {\rm e}^{-t\mathcal{L}}f\le 1$.
\end{enumerate}
Moreover, the inverse operator $\mathcal{L}^{-1}$ can be written  as
\[
\mathcal{L}^{-1}[f](x)=\int_\Omega {\mathbb K}(x,y)f(y)\dy\,,
\]
The kernel   $\mathbb{K}$ is called the Green function and we assume that there exist constants $\gamma\in (0,1]$ and $c_{0,\Omega},c_{1,\Omega}>0$ such that, for a.e. $x,y\in \Omega$\,:
\[\tag{K2}
c_{0}\,\p(x)\,\p(y) \le {\mathbb K}(x,y)\le \frac{c_{1}}{|x-y|^{N-2s}}
\left(\frac{\p(x)}{|x-y|^\gamma}\wedge 1\right)
\left(\frac{\p(y)}{|x-y|^\gamma}\wedge 1\right),
\]
where we adopt the notation $\delta(x):=\dist(x, \partial\Omega)$. We will also use $\phi_1$, the first eigenfunction of $\mathcal L$, and we know that $\phi_1\asymp \dist(\cdot,\partial\Omega)^\gamma$.  Further assumptions will be made in each statement, depending on the desired result we want, in particular (K4) that we have already mentioned.

\smallskip

\noindent $\bullet$ {\sc Sharp bounds.} Under these assumptions, we obtain existence and uniqueness of solutions with various properties, like time decay in $L^p$ spaces. We will not delve in this basic theory that is covered in the papers \cite{BSV2015, BV.Bdd2016}. We will stress here one of our main contributions in \cite{BFV2016}:  we  prove sharp upper and lower pointwise bounds for nonnegative solutions, both at the interior and close to the boundary. Indeed, we must pay close attention to the boundary behaviour, that turns out to be different for different operators in this class. However, only some options appear, as we describe next. Let us introduce first an important exponent
$$
\sigma=1\wedge\frac{2sm}{\gamma(m-1)}.
$$
Notice that $\sigma=1$ for the RFL and the CFL, but not always for the SFL unless $m=1$. The results that follow are taken from the last work, \cite{BFV2016}. In the next results (CDP) means the Cauchy-Dirichlet problem with zero lateral data, and solutions means dual weak solutions. We make the default assumptions (A1), (A2), and (K2) on $\mathcal{L}$.

\smallskip

\noindent\textbf{Case 1. Nonlocal operators with nondegenerate kernels.} We assume here moreover that the kernel of $\mathcal{L}$ is non degenerate at the boundary, namely
\begin{equation}\label{Operator.Hyp.lower.sigma}
\A f(x)=\int_{\RR^N}\big(f(x)-f(y)\big)K(x,y)\dy\,,\qquad\mbox{with } \inf_{x,y\in \Omega}K(x,y)\ge \underline{\kappa}_\Omega>0\,.
\end{equation}
Under these assumptions we can prove the following first version of the Global Harnack Principle.

\begin{theorem}\label{thm.GHP.PME.I}
Let (A1), (A2), (K2), and \eqref{Operator.Hyp.lower.sigma} hold.
Also, when $\sigma<1$, assume that $K(x,y)\le c_1|x-y|^{-(N+2s)}$ for a.e. $x,y\in \RR^N$
and that $\phi_1\in C^\gamma(\Omega)$.
Let $u\ge 0$ be a weak dual solution to the (CDP) corresponding to $u_0\in \rm L^1_{\phi_1}(\Omega)$. Then, there exist  constants $\underline{\kappa},\overline{\kappa}$, so that the following inequality holds:
\begin{equation}\label{thm.GHP.PME.I.Ineq}
\underline{\kappa}\, \left(\frac{t}{t+t_*}\right)^{\frac{m}{m-1}}\frac{\phi_1(x)^{\sigma/m}}{t^{\frac{1}{m-1}}}
\le \, u(t,x) \le \overline{\kappa}\, \frac{\phi_1(x)^{\sigma/m}}{t^{\frac1{m-1}}}
\end{equation}
for all $t>0$ and all $x\in \Omega$.
\end{theorem}

Here,  $t_*=\kappa_*\|u_0\|_{\rm L^1_{\phi_1}(\Omega)}^{-(m-1)}$, and this time will appear in the other theorems. For large times both lower and upper bounds are similar. We point out that the result holds in particular for the Restricted and Censored Fractional Laplacians, but not for the Spectral Fractional Laplacian. The lower bound is  false for $s=1$ (in view of the finite speed of propagation of the standard PME).

\medskip

\noindent\textbf{Case 2. Matching behaviour for large times.} We can prove that previous Global Harnack Principle for large times without using the non-degeneracy of the kernel, under the following conditions on $\sigma$ : either

(i) $\sigma=1$ (i.e., $2s>\gamma (m-1)/m$), or

(ii) $\sigma<1$, and we have an improved version of (K2)
\[\tag{K4}
{\mathbb K}(x,y)\asymp \frac{c_{1}}{|x-y|^{N-2s}}
\left(\frac{\p(x)}{|x-y|^\gamma}\wedge 1\right)
\left(\frac{\p(y)}{|x-y|^\gamma}\wedge 1\right),
\]
and the initial data are not small: $u_0\geq c\phi_1^{\sigma/m}$ for some $c>0$.

\begin{theorem}[Global Harnack Principle II]\label{thm.GHP.PME.II}
Let  (A1), (A2), and (K2) hold, and let $u\ge 0$ be a weak dual solution to the (CDP) corresponding to $u_0\in \rm L^1_{\phi_1}(\Omega)$. Assume that either (i) or (ii) above hold true. Then there exist  constants $\underline{\kappa},\overline{\kappa}>0$ such that the following inequality holds:
\begin{equation}\label{thm.GHP.PME.II.Ineq}
\underline{\kappa}\,  \frac{\phi_1(x)^{\sigma/m}}{t^{\frac{1}{m-1}}}
\le \, u(t,x) \le \overline{\kappa}\, \frac{\phi_1(x)^{\sigma/m}}{t^{\frac1{m-1}}}\qquad\mbox{for all $t\ge t_*$ and all $x\in \Omega$}\,.
\end{equation}
The constants $\underline{\kappa},\overline{\kappa}$ depend only on   $N,s,\gamma, m, \underline{\kappa}_0,\underline{\kappa}_\Omega$, and $\Omega$.
\end{theorem}

The conditions on $\sigma $ are sharp. Actually, the proof in the case $\sigma=1$  includes the classical PME (i.e., the non fractional equation, for which finite propagation holds, so that there can be no positive a priori lower bound for short times).

\medskip

\noindent\textbf{The case of a really degenerate kernel. } We assume moreover that
 we assume moreover that the kernel  of $\mathcal{L}$ exists and can be degenerate at the boundary (actually, excluding the local case, this is the most general assumption) in the form
\begin{equation}\label{Operator.Hyp.lower}
\A f(x)=P.V.\int_{\RR^N}\big(f(x)-f(y)\big)K(x,y)\dy\,,\quad\mbox{with $K(x,y)\ge c_0\phi_1(x)\phi_1(y)\; \,\,\forall\,x,y\in \Omega$} \,.
\end{equation}
This is an assumption that holds for the Spectral Fractional Laplacian operator. To our knowledge, precise information about the kernel of the SFL was not known before Lemma 3.1 of \cite{BFV2016}.

Note that, for small times, we cannot find matching powers for a global Harnack inequality (except for some special initial data), and such result is actually false for $s=1$ (in view of the finite speed of propagation of the PME). Hence, in the remaining cases, we have only the following general result.

\begin{theorem}[Global Harnack Principle III]\label{thm.GHP.PME.III}
Let (A1), (A2), (K2), and \eqref{Operator.Hyp.lower} hold.
Let $u\ge 0$ be a weak dual solution to the (CDP) corresponding to $u_0\in \rm L^1_{\phi_1}(\Omega)$. Then, there exist  constants $\underline{\kappa},\overline{\kappa}>0$, so that the following inequality holds:
\begin{equation}\label{thm.GHP.PME.III.Ineq}
\underline{\kappa}\,\left(\frac{t}{t+t_*} \right)^{\frac{m}{m-1}}\frac{\phi_1(x)}{t^{\frac{1}{m-1}}}
\le \, u(t,x) \le \overline{\kappa}\, \frac{\phi_1(x)^{\sigma/m}}{t^{\frac1{m-1}}}
\end{equation}
for all $t>0$ and all $x\in \Omega$\,.

\end{theorem}
This is what we call non-matching powers for the spatial profile at all times. The paper gives analytical and numerical evidence that such non matching behaviour does not happen in the associated elliptic problems, and came as a surprise to the authors. For some class of initial data, namely $u_0\le \varepsilon_0\phi_1$ we can prove that for small times
\[
\underline{\kappa}_0\,\left(\frac{t}{T}\right)^{\frac{m}{m-1}}\frac{\phi_1(x)}{t^{\frac{1}{m-1}}}
\le \, u(t,x) \le \overline{\kappa}_0\, T^{\frac{1}{m-1}}\frac{\phi_1(x)}{t^{\frac{1}{m-1}}}\qquad\mbox{for all $0\le t\le T$ and all $x\in \Omega$}\,.
\]

\noindent {\sc Numerics.} This work has been improved  in January 2017 with numerics done at BCAM Institute by my former student F. del Teso and collaborators, \cite{Cusi.numerics}, 2017, that validates the different behaviour types.

\medskip

\noindent $\bullet$ {\sc Asymptotic Behaviour. }
An important application of the Global Harnack inequalities of the previous section concerns the sharp asymptotic behavior of solutions. More precisely, we first show that for large times all solutions behave like the separate-variables solution $\mathcal{U}(t,x)={S(x)}\,{t^{-\frac{1}{m-1}}}$. The profile $S$ is the solution of an elliptic nonlocal problem. Then, whenever the Global Harnack Principle (GHP) holds, we can improve this result to an estimate in relative error.

\begin{theorem}[Asymptotic behavior]\label{Thm.Asympt.0}
Assume that $\A$ satisfies (A1), (A2), and (K2), and let $S$ be as above. Let $u$ be any weak dual solution to the (CDP).
Then, unless $u\equiv 0$,
\begin{equation}\label{Thm.Asympt.0.1}
\left\|t^{\frac{1}{m-1}}u(t,\cdot)- S\right\|_{{\rm L}^\infty(\Omega)}\xrightarrow{t\to\infty}0\,.
\end{equation}
\end{theorem}

We can exploit the (GHP) to get a stronger result, using the techniques of paper \cite{BSV2015}.

\begin{theorem}[Sharp asymptotic behavior]\label{Thm.Asympt}
Under the assumptions of Theorem \ref{Thm.Asympt.0}, assume that $u\not\equiv 0$.
Furthermore,
suppose that either the assumptions of Theorem \ref{thm.GHP.PME.I} or of Theorem \ref{thm.GHP.PME.II} hold.
Set $\mathcal{U}(t,x):=t^{-\frac{1}{m-1}}S(x)$.
Then there exists $c_0>0$ such that, for all $t\ge t_0:=c_0\|u_0\|_{{\rm L}^1_{\phi_1}(\Omega)}^{-(m-1)}$, we have
\begin{equation}\label{conv.rates.rel.err}
\left\|\frac{u(t,\cdot)}{\mathcal{U}(t,\cdot)}-1 \right\|_{{\rm L}^\infty(\Omega)} \le \frac{2}{m-1}\,\frac{t_0}{t_0+t}\,.
\end{equation}
We remark that the constant $c_0>0$ only depends on $N,s,\gamma, m, \underline{\kappa}_0,\underline{\kappa}_\Omega$, and $\Omega$.
\end{theorem}

\medskip

\noindent {\bf Comments on related work.} Construction of the solutions of the FPME on bounded domains with the SFL was already used in \cite{pqrv1, pqrv2} as an approximation to the problem in the whole space, but the regularity or asymptotic properties were not studied. Kim and Lee in \cite{KimLee11} study the Fast Diffusion range $m<1$ with a fractional Laplacian in a bounded domain and prove H\"older regularity and asymptotic behaviour. There is a developing literature on nonlocal nonlinear diffusion equations on domains.


\subsection{Fractional diffusion equations of $p$-Laplacian type}

We  report here about our work \cite{Vaz16}. It deals with a model of fractional diffusion involving a nonlocal version of the $p$-Laplacian operator, and the equation is
\begin{equation}\label{frplap}
\fbox{ $ \partial_t u+{\mathcal L}_{s,p} u=0, \displaystyle \qquad {\mathcal L}_{s,p}(u):=\int_{\ren}\frac{\Phi(u(y,t)-u(x,t))}{|x-y|^{N+sp}}\,dy=0 \ $}
\end{equation}
where $x\in \Omega\subset \ren$, $ N\ge 1,$ $\Phi(z)=c|z|^{p-2}z,$ $p\in (1,\infty)$ and $s\in (0,1)$. ${\mathcal L}_{s,p}$ is the Euler-Lagrange operator corresponding to a power-like functional with nonlocal kernel of the $s$-Laplacian type. The study of the equation is motivated by the recent increasing interest in nonlocal generalizations of the porous-medium equation. In the paper we cover the range $2<p<\infty$. Note that for $p=2$ we obtain the standard $s$-Laplacian heat equation, $u_t+(-\Delta)^{s} u=0$, which was discussed before; on the other hand, it is proved that in the limit $s\to1$ with $p\ne 2$, we get the well-known $p$-Laplacian evolution equation $\partial_t u=\Delta_p(u)$, after inserting a normalizing constant.

We  consider the equation in a bounded domain $\Omega\subset \ren$ with initial data
\begin{equation}\label{eq.ic.bdd}
u(x,0)=u_0(x), \qquad x\in \Omega,
\end{equation}
where $u_0$ is a nonnegative and integrable function. Moreover, we impose the homogeneous Dirichlet boundary condition that in the fractional Laplacian setting takes the form
\begin{equation}\label{eq.bc}
 u(x,t)=0 \qquad \mbox{ for all} \  x\in \ren, x\not\in \Omega,\quad \mbox{ and all $t>0$}.
\end{equation}
When then apply the integral operator  on the set of functions that vanish outside of $\Omega$

The first result of this paper  concerns the existence and uniqueness of a strong nonnegative solution to an initial-boundary value problem for \eqref{frplap} in  bounded domain $\Omega\subset\ren$, with zero Dirichlet data outside $\Omega\times(0,\infty)$. The boundedness of the solution is established after proving the existence of a special separating variable solution of the form
   $$
U(x,t)=t^{−1/(p-2)}F(x),
$$
called the friendly giant. The profile function $F(x)$ of the friendly giant solves the interesting nonlocal elliptic problem
$$
\int_{\ren}\frac{\Phi(F(y)-F(x))}{|x-y|^{N+sp}}dy=c \,F(x).
$$
 The friendly giant solution provides a universal upper bound and also gives the large-time behaviour for all the nonnegative solutions of initial-boundary value problems with homogeneous Dirichlet boundary conditions.

The fractional $p$-Laplacian has recently attracted the attention of many researchers  for its mathematical interest. See among other related works the papers by Caffarelli et al. \cite{CafChanVas}, Maz\'on et al. \cite{MRT}, Puhst \cite{Puhst}. Another approach in the form of non-local gradient dependent operators is taken by Bjorland-Caffarelli-Figalli in \cite{BjCaffFig}. The corresponding stationary equation is also studied in the literature, see previous references. Finally, the  work \cite{HyndL16} by Hynd and Lindgren deals with the doubly nonlinear model \ $|u_t|^{p-2} u_t+{\mathcal L}_{s,p} u=0$, that has a special homogeneous structure. Regularity and asymptotic behaviour follow.

\medskip


\section{\textsf{Further work on related topics}}
\label{sec.further}


\subsection{Diffusion with fractional time derivatives}

 Equations of the form
$$
D_t^\alpha u ={\mathcal L}u + f
$$
are another form of taking into account nonlocal effects. Here  $\mathcal L$ represents the diffusion process with long-distance effects in the family of the fractional Laplacian operators. The symbol $D_t^\alpha$ denotes the fractional time derivative. There are a number of variants of this concept, the most popular being maybe the {\sl  Caputo fractional derivative}, which was introduced by M. Caputo in  1967 \cite{Caputo} and reads
\begin{equation}
\fbox{ $
\displaystyle {}_{a}^{C}D_{t}^{\alpha }f(t)={\frac {1}{\Gamma (n-\alpha )}}
\int _{a}^{t}{\frac {f^{(n)}(\tau )\,d\tau }{(t-\tau )^{\alpha +1-n}}}, \quad n-1 <\alpha\le n$}
\end{equation}
Indeed, fractional time derivatives are the most elementary objects of Fractional Calculus, a branch of mathematical analysis that studies the possibility of taking real number powers (real number fractional powers or complex number powers) of the differentiation operator $    D = d / d x ,$ and the integration operator.

The foundations of the theory of fractional derivatives were laid down by Liouville in a paper from 1832. Different definitions use different kernels,     but all them make weighted averages in time.

\smallskip

Some recent work:  Dipierro and Valdinoci derive the linear time-fractional heat equation in 1D in a problem of neuronal transmission in cells, \cite{DPiVa17}; Allen, Caffarelli and Vasseur  study porous medium flow with both a fractional potential pressure and a fractional time derivative,\cite{ACV17}.

\subsection{Diffusion equations on Riemannian manifolds}

\noindent $\bullet$ The study of the heat equation posed on a Riemannian manifold comes from long time ago, the diffusive operator  being the Laplace-Beltrami operator
 \begin{equation}
\Delta_g(u)=\frac1{\sqrt{| g|}}\partial_i(g^{ij}\sqrt{| g|} \partial_j u)\,,
 \end{equation}
with the usual Riemannian notations. Thus, the heat kernel was studied in \cite{Yau78, Dav87, Grig87}, random walks and Brownian motion on manifolds are studied in \cite{Var87}. There are lots of recent works, an example is \cite{Grig06}.

\noindent $\bullet$ Generalization of the Caffarelli-Silvestre extension method allows to define extensions and boundary operators of the fractional Laplacian type when $ M$ is the boundary  of a conformally compact Einstein manifold. Combining geometrical and PDE approaches,  Chang and Gonz\'alez in \cite{ChGonz2011}  related the original  definition of the conformal fractional Laplacian coming from scattering theory to a Dirichlet-to-Neumann operator for a related elliptic extension problem for $M$, see also \cite{Gonzsurv}. It is possible then to formulate fractional Yamabe-type problems for the conformally covariant operators $P_\gamma$, \cite{GonzQing}. For more recent work in this problem see \cite{JX2014, DSV2017}.

\noindent $\bullet$ Work on the Porous Medium Equation on manifolds was done in the 2000s, like \cite{BG05}, \cite{Vazpme}; fast diffusion was treated in \cite{BGV10}; general Aronson-B\'enilan estimates and entropy formulae for porous medium and fast diffusion equations on manifolds were obtained in \cite{LuNiVaVi}.

\noindent $\bullet$ Recently the author studied the PME on the hyperbolic space \cite{Vazhyp15}, and constructed the fundamental solution and proved the asymptotic convergence and free boundary propagation rates. The fact that the fundamental is not explicit or self-similar is not pleasant, and looking for some higher symmetry properties a remarkable object appeared to play an asymptotic role. Namely, there exists  nonnegative weak solution  of the PME defined on the whole of $\mathbb{H}^N$ for all $t>0$, that has a strong algebraic  structure. In the Poincar\'e upper half-space representation it is given by the formula
\begin{equation}\label{gtw.mn}
U(x,y,t)^{m-1}=a\frac{(\log(ct^{\gamma}y))_+}{t}
\end{equation}
with $m,N>1$, $1/a=m(N-1)$ and $1/\gamma=(N-1)(m-1)$, and $x\in \re^{N-1}$, $y>0$.  Note that $U$ has zero initial trace at $t=0$ on the half-space, but it has a singularity as trace at $y=+\infty$, which corresponds to a singularity at the North Pole in the standard Poincar\'e ball representation. Therefore, we can say that the special solution (or geometrical soliton) $U$ comes from the infinite horizon and expands to gradually to fill the whole space; for any $t>0$ it has a support limited by a family of horospheres $\{(x,y): y=(1/c)t^{-\gamma}\}$. Recall that geodesic distance is given by the formula $ds^2=dy^2/y^2$. A detail for analysts: $U$ represents an example of non-uniqueness of nonnegative solutions for the Cauchy Problem in hyperbolic space.

This family of special solutions is the pattern to which all solutions with compactly supported initial data are proved to converge as $t\to\infty$. Accordingly, we get the following sharp estimates for all solutions $u\ge 0$ with compact support:
$$
\|u(\cdot,t)\|_\infty\sim ct^{-1/(m-1)}\log(t)^{1/(m-1)}, \quad S(t)\sim \gamma\log(t)
$$
where $S(t)$ is the location of the free boundary measured in geodesic distance.

\noindent $\bullet$ The asymptotic analysis of PME flows is extended to more general manifolds in \cite{GMV}. The work on the asymptotic behaviour on hyperbolic space  is extended to fast diffusion by Grillo et al. in \cite{GMP1}. In another direction, Amal and Elliott \cite{AmEll6} study fractional porous medium equations  on evolving surfaces, a very novel subject.


\subsection{Diffusion in inhomogeneous media}

We have already mentioned the inhomogeneity of the medium as a reason for the introduction of coefficients in the passage from the heat equation to the parabolic class. In view of the important practical consequences, there is no surprise in finding coefficients appear in most of the models we have considered above, both linear and nonlinear, local and nonlocal. Let us just mention
some well-known references like \cite{KR81, ReV09, KRV10}, where it appears as diffusion with weights, see  also the more recent \cite{GMPw} with two weights. There is an interesting connection between weighted diffusion in Euclidean space and Laplace-Beltrami diffusion on manifolds, that has been studied in \cite{Vazhyp15} and is being further investigated.


\subsection{Drift diffusion equations with or without fractional terms}

The  main equation in this case is an equation for a scalar unknown $\theta$ driven by the equation
\begin{equation}
\partial_t \theta + {\bf v}\cdot \nabla_x \theta={\mathcal L}(\theta)\,,
\end{equation}
where ${\bf v}\cdot \nabla_xT$ is the convective term with velocity vector ${\bf v}$, and ${\mathcal L}(\theta)$ is the diffusion operator; that diffusion can be linear of nonlinear, local or nonlocal. A very important aspect is the relation of ${\bf v}$ to the rest of the variables. Thus, when ${\bf v}={\bf v}(x,t)$ is a given function, no essential new problems arise if ${\bf v}$ is smooth. But serious difficulties happen for nonsmooth ${\bf v}$.  All this is reflected for instance in the seminal paper by Caffarelli and Vasseur \cite{CaffVass} where the equation is
$$
\partial_t \theta + {\bf v}\cdot \nabla_x \theta+(-\Delta)^{1/2}(\theta)=0
$$
and  ${\bf v}$ is a divergence-free vector field. In the popular quasi-geostrophic model ${\bf v}$  is given in terms of $\theta $ that makes the problem more involved, but the results stated in \cite{CaffVass}  do not depend upon such dependence. The proof of regularity needs to establish delicate local energy estimates, despite the fact that the diffusion operator $(-\Delta)^{1/2}$ is non-local. It also uses DeGiorgi's methods in an essential way. There is much work in geostrophic flows, like \cite{MT, KNS1, KNV}.

There is a vast literature on this important issue. In some cases ${\bf v}$ is given by Darcy's law in an incompressible fluid, and the papers refer to the problem as ``flow in porous media'', like \cite{CCGO}. Let us point out that this use is quite different from our use of ``porous media'' in the present paper, the difference being often stressed by calling their use ``incompressible flow in porous media'', \cite{Fabrie86}.

\subsection{Other}

-Minimal surfaces are an important subject which uses many methods of the nonlinear elliptic and parabolic theory, Recently, it has developed a new branch, nonlocal minimal surfaces. Work on both aspects is reported in detail in another contribution to this volume by Cozzi and Figalli
\cite{CF-LN}. Related items are fractional perimeters and fractional phase transition interfaces. We will not enter into that area.

-The study of the combined effects of diffusion and aggregation is a very active field where the methods of diffusion in its different forms must be combined with the counter mechanism of attraction. We refer to the contribution by Calvez, Carrillo and Hoffmann to this volume, \cite{CCH}.

- In the study of nonlinear diffusion  we have chosen to present almost exclusively equations with diffusion terms in divergence form. There is a large body of work involving Fully Nonlinear Parabolic Equations (they are non-divergence equations), both elliptic and parabolic. We will not touch such theories here.

\section*{Addendum and final comment}

Here is the complete Wikipedia list of diffusion topics:

Anisotropic diffusion, also known as the Perona-Malik equation, enhances high gradients; Anomalous diffusion, in porous medium; Atomic diffusion, in solids;  Brownian motion, for example of a single particle in a solvent;  Collective diffusion, the diffusion of a large number of (possibly interacting) particles;  Eddy diffusion, in coarse-grained description of turbulent flow;
Effusion of a gas through small holes;  Electronic diffusion, resulting in electric current;  Facilitated diffusion, present in some organisms; Gaseous diffusion, used for isotope separation; Heat flow, diffusion of thermal energy; It\=o-diffusion, continuous stochastic processes; Knudsen diffusion of gas in long pores with frequent wall collisions;  Momentum diffusion, ex. the diffusion of the hydrodynamic velocity field;   Osmosis is the diffusion of water through a cell membrane; Photon diffusion;  Random walk model for diffusion; Reverse diffusion,  against the concentration gradient, in phase separation; Self-diffusion;  Surface diffusion, diffusion of adparticles on a surface;
Turbulent diffusion, transport of mass, heat, or momentum within a turbulent fluid.

\medskip

\noindent $\bullet$ The reader may wonder whether mathematical diffusion is a branch of applied mathematics? In principle it would seem that the answer is an obvious yes, and yet it is not so clear. As we have seen, the mathematical theories of diffusion have developed into a core knowledge in pure mathematics, that encompasses several branches, from analysis and PDEs to probability, geometry and beyond. We hope that the preceding pages will have convinced the reader of this trend, and will also motivate him/her to pursue some of many avenues open towards the future.

\medskip

{\sc Acknowledgments. } This work was partially supported by Spanish Project MTM2014-52240-P. The text is based on series of lectures given at the CIME Summer School held in Cetraro, Italy, in July 2016. The author is grateful to the CIME foundation  for the excellent organization. The author is also very grateful to his collaborators mentioned in the text for an effort of many years. Special thanks are due to F. del Teso, N. Simonov and D. Stan for a careful reading and comments on the text.

\

\newpage


{\small
\bibliographystyle{amsplain}

}

\medskip

\medskip

{
\noindent {\sc Address}: Juan Luis V\'azquez\\[2pt]
   {Department de Matem\'aticas}\\[2pt]
   {Universitas Aut\'onoma de Madrid}\\[2pt]
   {28049 Madrid,  Spain}\\[2pt]
   {\sl e-mail:} {juanluis.vazquez@uam.es}
}

\vskip .5cm

2010 \textit{\sc  Mathematics Subject Classification.}
26A33, 
35K55, 
35K65, 
35S10. 

\medskip

\textit{\sc Keywords and phrases.} Diffusion, nonlinear diffusion, nonlocal diffusion, fractional operators.

\medskip

\sl Figures obtained from Google and personal sources.\rm

\end{document}